\documentclass[11pt]{article}
\usepackage{amsmath, amssymb, eucal, latexsym}
\usepackage[cp1251]{inputenc}
\usepackage[russian]{babel}

\usepackage{graphicx,a4,color}%

\newtheorem{theorem}{Theorem}

\newtheorem{lemma}{Lemma}

\newtheorem{remark}{Remark}

\begin{document}

\begin{center}

\LARGE{\bf Discrimination between close hypotheses about Weibull and log-Weibull type distributions by the higher order statistics}
\vspace{5mm}

\Large{Igor Rodionov}
\vspace{3mm}

\normalsize{Moscow State University, Faculty of Mathematics and Mechanics}

{Moscow Institute of Physics and Technics, Faculty of Innovations and High Technologies}

{\small e-mail: vecsell@gmail.com}

\end{center}

\section{Introduction.}

The problem of discrimination between similar tailed distributions appears
in many applications of extreme value theory, for example, linked with high
 risk insurance problems, see Gupta and Kundu (2003), Kundu and Raqab (2007).
  In addition, distributions of the intermediate values are often well-behaved
  for modeling by standard distributions that differ from asymptotical behavior of the tails.
   It seems that theory of contiguity of probability measures, see Roussas (1972), is an important instrument of discriminating between families of distributions with close tails and estimating of power of various discriminating tests. In presented work asymptotical behavior of the ratio of likelihoods

\begin{equation}R_n(u) = \frac{L(X_{n,n}, \ldots, X_{n-k_n+1,n};
\gamma+t(k_n, u))}{L(X_{n,n}, \ldots, X_{n-k_n+1,n}; \gamma)}\label{r}\end{equation} as
$n\rightarrow\infty, k_n\rightarrow\infty,
\frac{n}{k_n}\rightarrow\infty$ is considered for the family of distributions with the infinite right point

\begin{equation}f(x,\gamma) = \exp(-S(x,\gamma)),\label{f}\end{equation}
 where the function $S(x, \gamma)$ is positive starting from some $x_0(\gamma)>0.$ Suppose the function $S(x, \gamma)$ is strictly monotone and 4 times continuously differentiable by $x$ for all $\gamma$ and $x > x_0(\gamma) > 0$. There are two types of regularity conditions that will be imposed on the function $S(x, \gamma)$ in our work. The first type makes possible to find the asymptotical behavior of $R_n(u)$ as $n\to \infty$ for the family of probability densities $f_1(x, \gamma)$ with $S(x,\gamma)\gtrsim x^\varepsilon,$ $\varepsilon>0.$ Now introduce some examples of the families of probability densities that satisfy the first type of regularity conditions. Firstly, it should be mentioned the family of Weibull-type probability densities \[f_W(x,\gamma) = C(\gamma)\exp( - x^\gamma), \;x\geq 0,\;\gamma>1.\] One another example is the family of normal densities, where the parameter $\gamma$ is the variance of concerned distribution: $$f_{N}(x,\gamma) = \frac{1}{\sqrt{2\pi\gamma}} e^{-\frac{x^2}{2\gamma}},\;\gamma>0.$$ The family of Gumbel-type probability densities also satisfies regularity conditions: \[f_{G}(x, \gamma) = \gamma e^{\gamma x} e^{- e^{\gamma x}},\;x \geq 0, \;\gamma>0.\]
 If we want to expand our result to a wider family of probability densities $f_2(x, \gamma)$ with $S(x, \gamma)\gtrsim (\ln x)^\varepsilon,\ \varepsilon>1,$ then we impose the second type of regularity conditions on the concerned family of densities, that are slightly stricter than the first one. Note that the probability density $f(x, \gamma)$ with $S(x, \gamma) \sim \ln x$ belongs to the Frech\'{e}t maximum domain of attraction as opposed to the families of probability densities with $S(x, \gamma)\gtrsim (\ln x)^\varepsilon,\ \varepsilon>1,$ belonging to the Gumbel maximum domain of attraction (this fact is proved later). So, the second type of regularity conditions makes possible to consider practically all distributions with finite right point belonging to the Gumbel domain of attraction. Provide some examples of the families of probability densities that satisfy the second type of regularity conditions. The first important example is the family of log-Weibull-type probability densities: \[f_{LW}(x,\gamma) = C(\gamma)\exp( - (\ln x)^\gamma), \;x\geq 1,\;\gamma>1.\] The family of Weibull-type probability densities satisfy the regulrity conditions as before: \[f_{W}(x,\gamma) = C(\gamma)\exp( - x^\gamma), \;x\geq 0,\;\gamma>1.\] Another example is the family of the log-normal densities: \[f_{LN}(x,\gamma) = \frac{1}{\sqrt{2\pi\gamma}x} e^{-\frac{(\ln x)^2}{2\gamma}},\;\gamma>0, x\geq 0.\]
The method of the likelihood ratio as well as the method of the ratio of maximized likelihoods (RML-test) is well-known and often applied for discriminating between close types of distributions. In this connection the following works (Antle and Bain (1969), Antle and Dumonceaux (1973),
Dumonceaux at el. (1973)) should be mentioned. The method of the ratio of maximized likelihoods is applied for discriminating between Weibull-type distributions in the works Gupta and Kundu (2003), Gupta at el. (2001), Dey and Kundu (2012) and other.
In addition the wide families of distributions belonging to the Gumbel maximum domain of attraction are discriminated in the works
de Haan at el. (2009) and Girard at el. (2011).\\
\\
Consider the method of the likelihood ratio applied to $k_n$ higher order statistics. A similar method is applied in the work Dey and Kundu (2012), but there the lower order statistics for distributions with finite left point are considered, whereas the higher order statistics are significant for extreme value theory in this case. In this work the asymptotical behavior of the likelihood ratio for the family of probability densities $f(x, \gamma)$ satisfying some regularity conditions is considered.

\section{Main results.}
Return to the mentioned above regularity conditions, that are imposed on the function $S(x, \gamma)$.
 Since we consider the close hypotheses with respect to the parameter $\gamma,$
 it is natural to make propositions about the general form of the function $S(x, \gamma):$
 for example, we can consider functions like $\gamma S(x)$ or $(S(x))^\gamma.$
  And it is unreasonably to deal with functions like $f(x) + g(x,\gamma),$
  since the summand that is not dependent on $\gamma$ is cancelled in the ratio of likelihoods, therefore we will consider in the sequel
  that $S(x,\gamma)$ is equal to $g(x,  \gamma)$ only.
So it is reasonably to suppose about behavior of the partial derivative of $S(x, \gamma)$ by $\gamma$ as $x\to\infty$ the following: $\lim\limits_{x\rightarrow +\infty}\frac{\ln \left|\frac{\partial S(x,\gamma)}{\partial \gamma}\right|}{\ln S(x,\gamma)}= 1.$ \\
\\Impose the following regularity conditions on the function $S(x,\gamma)$:
\begin{enumerate}
\item[\bf A1] There exist $\epsilon = \epsilon(\gamma)>0$ for all $\gamma$ that $\frac{S(x, \gamma)}{x^{1+\epsilon}}\rightarrow +\infty$ as $x\rightarrow +\infty.$

\item[\bf A2] All partial derivatives of the function $S(x, \gamma)$ up to the fourth derivative are not equal to 0 or equal to 0 identically for all $\gamma$ as $x>x_1(\gamma)>x_0(\gamma).$ In addition all partial derivatives of the function $S(x,\gamma)$ up to the third derivative have a finite or an infinite limit as $x\rightarrow +\infty.$

\item[\bf A3] Suppose $T(x, \gamma)$ is the arbitrary partial derivative of the function $S(x, \gamma)$ up to the third derivative. There holds for all $\gamma$ and $k=1, 2, 3$: $\lim\limits_{x\rightarrow +\infty}\frac{\ln \left|\frac{\partial^{k} S(x,\gamma)}{\partial \gamma^k}\right|}{\ln S(x,\gamma)}= 1,$ if these partial derivatives are not equal to 0 identically. All expressions $\frac{\frac{\partial \ln |T(x, \gamma)|}{\partial x}} {\frac{\partial \ln S(x, \gamma)}{\partial x}}$ have a limit for all $\gamma$ as $x\to +\infty,$ if $T(x, \gamma)$ is not equal to 0 identically.
\end{enumerate}
 Note that it follows from the first regularity condition that $\lim\limits_{x\rightarrow +\infty}\frac{\ln \frac{\partial S(x, \gamma)}{\partial x}}{\ln S(x, \gamma)}>0.$ Denote
\[a_{n/k_n} =
\overline{F}^{\leftarrow}\left(\frac{k_n}{n}, \gamma\right),\] where $F(x,\gamma)$ is the distribution function with the density $f(x, \gamma),$ $\overline{F}(x, \gamma) = 1 - F(x, \gamma),$ and $\overline{F}^{\leftarrow}(x, \gamma) = \inf\{t:\; \overline{F}(x, \gamma) = t\}.$ Denote $N(a, b)$ as the normal random variable with the mathematical expectation $a$ and the variance $b.$ Notation $\xi_n \xrightarrow{d} \eta$ $(\xi_n \xrightarrow{P} \eta)$ means that the random sequence $\{\xi_n\}$ converges in distribution (in probability) to the random variable $\eta$ as $n\to +\infty.$ Also we write $S^{''}_{x\gamma}(a_{n/k_n},\gamma + t(k_n, u))$ instead of $\left.\frac{\partial^2 S(x, \gamma)}{\partial x \partial \gamma}\right|_{(x, \gamma) = (a_{n/k_n}, \gamma+t(k_n, u))}$ and for all other partial derivatives of the function $\widetilde{S}(x, \gamma)$ analogously. In the following theorem limit distribution of the likelihood ratio $R_n(u)$ on the assumption of slow increase of the sequence $k_n$ is found:

\begin{theorem} Suppose the family of probability densities $\{f(x, \gamma)\}$ satisfies regularity conditions A1--A3, $n\rightarrow\infty,$ $k_n\rightarrow\infty,$ and there exists $\varepsilon,$ $0<\varepsilon<2,$ that
\begin{equation}\label{cond}\lim\limits_{n\to\infty}\frac{k_n}
{(\ln\frac{n}{k_n})^{\varepsilon}} = 1.\end{equation}
Suppose also
$t(k_n, u):=\frac{u}{\sqrt{k_n}}\frac{S'_{x}(a_{n/k_n},\gamma)}{S^{''}_{x\gamma}(a_{n/k_n},\gamma)},$ then
\[\ln R_n(u)\xrightarrow{d} N\left(-\frac{u^2}{2},
u^2\right) \text { as } n\to\infty\]
\end{theorem}
The parameter $u$ in notation $t(k_n, u)$ will be skipped in the sequel.\\
\\
\begin{remark} Note that if we suppose $k_n = n,$ i.e. the concerned likelihoods would depend on the whole sample, and take $t(k_n) =
\frac{u}{\sqrt{n}},$ then we come to a classical result of contiguity theory (see theorem 2.4.3 in Roussas (1972)).
\end{remark}
Since the proof of the theorem 1 seems technically complicated,
let's provide the scheme of the proof. Lemma 1
is used to express the ratio of likelihoods $R_n(u)$ given $X_{(n-k_n)} = q$ in explicit form.
Since we try to find conditional distribution of the likelihood ratio $R_n(u)$ given $X_{(n-k_n)},$ then
we prove the asymptotical properties of the $(n-k_n)$th order statistic $X_{(n-k_n)}$ and the quantile
$a_{n/k_n}$ in Lemma 3. Then the central limit theorem (\ref{N01}) is applied for independent
given $X_{(n-k_n)} = q$ and identically distributed random variables $Y_i$
emerging in the expression of the logarithm of the likelihood ratio. 
Using Lemma 2 we find the mathematical expectation (\ref{Ey1}) 
and the variance (\ref{Dy1})
 of the random variable $Y_1$. 
And in the end the residual terms in the logarithm of the likelihood ratio (\ref{123}) are taken into account that leads to a final answer.\\
\\
Now consider the second type of regularity conditions imposing on the functions $S(x, \gamma).$ Impose the following regularity conditions on the function $S(x,\gamma)$:
\begin{enumerate}
\item [\bf B1] There exists $\epsilon = \epsilon(\gamma)>0$ for all $\gamma$ that $\frac{S(x, \gamma)}{(\ln x)^{1+\epsilon}}\rightarrow +\infty$ as $x\rightarrow +\infty.$

\item [\bf B2] There exists $\delta = \delta(\gamma),$ $0\leq \delta <1,$ for all $\gamma$ that $\lim\limits_{x\rightarrow +\infty}\frac{\ln \frac{\partial S(x,\gamma)}{\partial x}}{(S(x,\gamma))^{1-\delta}} = 0.$

\item [\bf B3] All partial derivatives of the function $S(x, \gamma)$ up to the fourth derivative are not equal to 0 or equal to 0 identically for all $\gamma$ as $x>x_1(\gamma)>x_0(\gamma).$ In addition all partial derivatives of the function $S(x,\gamma)$ up to the third derivative have a finite or an infinite limit as $x\rightarrow +\infty.$

\item [\bf B4] Suppose $T(x, \gamma)$ is the arbitrary partial derivative of the function $S(x, \gamma)$ up to the third derivative. There holds for all $\gamma$ and $k=1, 2, 3$: $\lim\limits_{x\rightarrow +\infty}\frac{\ln \left|\frac{\partial^{k} S(x,\gamma)}{\partial \gamma^k}\right|}{\ln S(x,\gamma)}= 1,$ if these partial derivatives are not equal to 0 identically. All expressions $\frac{\frac{\partial \ln |T(x, \gamma)|}{\partial x}} {\frac{\partial \ln S(x, \gamma)}{\partial x}}$ have a limit for all $\gamma$ as $x\to +\infty$ and $\lim\limits_{x\to+\infty}\frac{\ln |T(x, \gamma)|} {S(x, \gamma)} = 0$ for all $\gamma,$ if $T(x, \gamma)$ is not equal to 0 identically.
\end{enumerate}
Denote
\[H(x) = \sqrt{k_n}t(k_n)\left(\frac{\int\limits_{x}^{\infty}S_{\gamma}(y,\gamma) \exp(-S(y,\gamma)) dx}{ \int\limits_{x}^{\infty} \exp(-S(y,\gamma))dx} - S_{\gamma}(x,\gamma) - \right.\]

\[\left. - \frac{S_{x\gamma}(x,\gamma)}{S_x(x,\gamma)} - \frac{S_{xx\gamma}(x, \gamma)}{(S_x(x, \gamma))^2} + 2 \frac{S_{xx}(x, \gamma)S_{x\gamma}(x, \gamma)}{(S_x(x, \gamma))^3}\right).\]
\medskip
It follows from the statement of Theorem 2 and its proof, that $H(x) =  - \frac{S_{xx}(x, \gamma)}{S^2_x(x, \gamma)}(1+o(1))$ as $x\to+\infty.$

\begin{theorem}

 Suppose the family of probability densities $f(x, \gamma)$ satisfies the regularity conditions B1--B4, $n\rightarrow\infty,$ $k_n\rightarrow\infty,$ and there exists $\varepsilon,$ $0<\varepsilon<2,$ that
\begin{equation}\label{cond2}\lim\limits_{n\to\infty}\frac{k_n}
{n^{\varepsilon}} = 1.\end{equation}
Suppose also
$t(k_n)=\frac{u}{\sqrt{k_n}}\frac{S'_{x}(a_{n/k_n},\gamma)}{S^{''}_{x\gamma}(a_{n/k_n},\gamma)}$ then
\[ \ln R_n(u) - \sqrt{k_n}H(a_{n/k_n}) \xrightarrow{d}  N\left(-\frac{u^2}{2},
u^2\right)\]
as $n\to +\infty.$
\end{theorem}

\begin{remark}
It should be noted, that the asymptotical behavior of the likelihood ratio $R_n(u)$ in  Theorem 1 and Theorem 2 does not depend on $\gamma,$ so it is convenient for construction of the criterion, see Chibisov (2009).
\end{remark}

\section{Auxiliary lemmas}
 In this section we state several auxiliary lemmas that are used in the proofs of Theorem 1 and Theorem 2. Lemma 1 makes possible to express the likelihood ratio in explicit form such as a function of $(n-k_n)$th order statistic.

\begin{lemma}(lemma 3.4.1 in de Haan (2006)) Let $X, X_1, \ldots, X_n$ be independent and identically distributed random variables with common distribution function $F(x)$, and let $X_{(1)}\leq \ldots \leq X_{(n)}$ be the $n$th order statistics. The joint distribution of $\{X_{(i)}\}_{i=n-k+1}^{n}$ given $X_{(n-k)}=q$
for some $k = 1 \ldots n-1$ equals the joint distribution of the set of order statistics
$\{X_{(i)}^{*}\}_{i=1}^{k}$ of independent and identically distributed random variables
$\{X_i^*\}_{i=1}^{k}$ with common distribution function
\[F_q(x) = P(X \leq x| X>t) = \frac{F(x) - F(q)}{1 - F(q)}, \ \ \
x>q.\]
\end{lemma}
 The following lemma (Rodionov (2014)) is extension of the Laplace method of estimating integrals and makes possible to find the asymptotic form of integrals emerging in presented work (see, for example, Fedoruk (1977)).
\begin{lemma} Consider the behavior of the integral $F(q) =
\int\limits_{q}^{+\infty}\exp(-S(x))dx$ as $q\to \infty.$ Let $S(x)$ satisfy the following regularity conditions:
\begin{enumerate}
\item[\bf C1] $\frac{S(x)}{\ln x}\rightarrow +\infty$ as $x\rightarrow +\infty.$

\item[\bf C2] $S(x)$ is strictly monotone and three times continuously differentiable starting from some $x_0>0.$
In addition, the first, the second and the third derivatives of the function $S(x)$ have a finite or an infinite limit as $x\rightarrow +\infty.$

\item [\bf C3] For $k=1, 2, 3$ $\lim\limits_{x\rightarrow +\infty} \frac{\ln \left|\frac{d^{k} S(x)}{dx^k}\right|}{S(x)} = 0.$
\end{enumerate}
Then  \[F(q)=
\exp(-S(q))\left(\sum\limits_{k=0}^{2} c_k +
o(c_2)\right) \ \ \ as\ \ q\to+\infty\]
where $c_k = M^k\left.\left(\frac{1}{S'(x)}\right)\right|_{x=q}$ and $M
= \frac{1}{S'(x)}\frac{d}{dx}$.
\end{lemma}
It is not difficult to prove using this lemma that the class of distributions satisfying regularity conditions A1--A3 (and B1--B4) belongs to the Gumbel maximum domain of attraction. The distribution function of such distributions may be represented in form of Von Mises (see Embrechts at el. (1997)):

\begin{equation}1 - F(x) = d(x) \exp\left( - \int\limits_{x'}^{x} \frac{g(t)}{a(t)}dt\right),\label{VonMises}\end{equation}
where $x' \geq 0,$ $a(x)$ is the positive and absolutely continuous function on $[x'; \infty),$ $a'(x) \to 0,$ $d(x) \to c > 0,$ $g(x) \to 1$ as $x \to x_F,$ and $x_F$ is the right point of concerned distribution. On the other hand, one can derive for the concerned family of probability densities $f(x, \gamma) = \exp(-S(x,\gamma))$  (since we differentiate only in respect to $x,$ we omit the parameter $\gamma$ here):

$$ 1 - F(x) =\exp(-S(x)) \frac{1}{S'(x)} \left(1 - \frac{S''(x) }{(S'(x))^2} + o\left(\frac{S''(x)}{(S'(x))^2}\right)\right).$$
It is proved in Rodionov (2014), that under the conditions of Theorem 1 and Theorem 2 $\lim\limits_{x\to\infty}\frac{S''(x) }{(S'(x))^2} = 0.$ It follows from the L'Hospital rule, that \[\lim\limits_{x\to\infty} \frac{S'(x)}{1/x} = \lim\limits_{x\to\infty} \frac{S(x)}{\ln x} = +\infty.\] So we can write
\[d(x) =\left(1 - \frac{S''(x) }{(S'(x))^2} + o\left(\frac{S''(x)}{(S'(x))^2}\right)\right)\]
for the function $d(x)$ in (\ref{VonMises}).
So, it should be proved only, that there exist $x',$ $a(x)$ and $g(x)$ in (5) that $S(x) + \ln S'(x) = \int\limits_{x'}^{x} \frac{g(t)}{a(t)}dt.$ Differentiate this equation in respect to $x: $
$$S'(x) + \frac{S''(x)}{S'(x)} = S'(x)\left(1 + \frac{S''(x)}{(S'(x))^2}\right) = \frac{g(x)}{a(x)}.$$
 It is mentioned earlier, that $\lim\limits_{x\to\infty}\frac{S''(x) }{(S'(x))^2} = 0$ under regularity conditions A1--A3 and B1--B4, so let $g(x) = 1 + \frac{S''(x)}{(S'(x))^2}.$ Since $\left(\frac{1}{S'(x)}\right)'= - \frac{S''(x)}{(S'(x))^2},$ let $a(x)$ be equal $\frac{1}{S'(x)},$ it is the positive function starting from some $x_0>0$ by the condition C2. Thus we state, that there exist $a(x)$ and $g(x)$ such that $S(x) + \ln S'(x) = \int\limits_{x'}^{x} \frac{g(t)}{a(t)}dt + C_1,$ where $C_1$ is a certain constant. To finish the proof of the state, that the family of concerned distributions belongs to the Gumbel maximum domain of attraction, take $x' = x_0$ and move the constant $e^{C_1}$ to $d(x).$\\
 \\
In the following lemma asymptotic distribution of the $(n-k_n)$th order statistic $X_{(n-k_n)}$ is considered.
\begin{lemma}There holds under the conditions of Theorem 1

$$X_{(n-k_n)} - a_{n/k_n}\xrightarrow{P} 0.$$
\end{lemma}
\textbf{Proof.} According to theorem 2.2.1 from de Haan (2006),
\begin{equation}\label{lemm}\sqrt{k_n}\frac{X_{(n-k_n)} -
U\left(\frac{n}{k_n}\right)}{\frac{n}{k_n}
U'\left(\frac{n}{k_n}\right)} \xrightarrow{d}
N(0,1),\end{equation} where
$U=\left(\frac{1}{1-F}\right)^{\leftarrow}.$ According to Lemma 2,
$$a_{n/k_n} =
\arg\left\{t: \;\int\limits_{t}^{\infty}\exp(-S(x,\gamma))dx =
\frac{k_n}{n}\right\} \xrightarrow[n\rightarrow\infty]{}$$
\begin{equation}\label{ank}\xrightarrow[n\rightarrow\infty]{} \arg
\left\{t:\; \frac{\exp(-S(t,\gamma))}{S'(t,\gamma)} =
\frac{k_n}{n}\right\}.\end{equation}
In addition,
$U(\frac{n}{k_n}) = a_{n/k_n}.$ Then $U'(t) =
\frac{[1-F(U(t))]^2}{F'(U(t))}$ and $U'(\frac{n}{k_n}) =
\frac{(k_n/n)^2}{F'(a_{n/k_n})}.$ But $F'(a_{n/k_n}) =
\exp(-S(t,\gamma)),$ where $t = a_{n/k_n}$ such that
$\frac{\exp(-S(t,\gamma))}{S'(t,\gamma)} = \frac{k_n}{n}.$ It follows from here that
 $F'(a_{n/k_n}) = \frac{k_n}{n}S'(a_{n/k_n},\gamma).$
Eventually obtain

\begin{equation}\frac{\sqrt{k_n}(X_{(n-k_n)} - a_{n/k_n})}{\frac{1}{S'(a_{n/k_n},\gamma)}} \xrightarrow {d} N(0,1).\label{N}\end{equation}
It follows from the regularity conditions A1 and A2 that $\lim\limits_{n\to +\infty} S'(a_{n/k_n}, \gamma) = +\infty,$ so we obtain using (\ref{N}) \[X_{(n-k_n)} - a_{n/k_n}\xrightarrow{P} 0.\]
The proof of Lemma 3 is completed.

\section{Proofs}
{\bf The proof of Theorem 1.}\\
\\
Firstly write the likelihood $L(X_{(n)}, \ldots,
X_{(n-k_n+1)}; \gamma)$ using Lemma 1

\[L(X_{(n)}, \ldots, X_{(n-k_n+1)}; \gamma) =
\frac{\prod\limits_{i=0}^{k_n-1} \exp(-S(X_{(n-i)},\gamma))}{\left (
\int\limits_{X_{(n-k_n)}}^{+\infty} \exp(-S(x,\gamma)) dx
\right )^{k_n}}.\]
So the ratio of likelihoods is the following:

$$R_n(u) = \frac{\prod\limits_{i=0}^{k_n-1} \exp[-S(X_{(n-i)},\gamma + t(k_n)) +
S(X_{(n-i)},\gamma)]} {\left ( \int\limits_{X_{(n-k_n)}}^{+\infty}
\exp[-S(x,\gamma)] dx\right)^{-k_n}\cdot
\left(\int\limits_{X_{(n-k_n)}}^{+\infty}
\exp[-S(x,\gamma+t(k_n))]dx \right)^{k_n}}.$$
Using Lemma 2, obtain

$$\int\limits_{X_{(n-k_n)}}^{+\infty}
\exp[-S(x,\gamma)] dx = $$
$$ = \frac{\exp(-S(X_{(n-k_n)},\gamma))}{S'_x(X_{(n-k_n)},\gamma)} \left(1 - \frac{S_{xx}^{''}(X_{(n-k_n)},\gamma) }{(S^{'}_x(X_{(n-k_n)},\gamma))^2} + o\left(\frac{S_{xx}^{''}(X_{(n-k_n)},\gamma)}{(S^{'}_x(X_{(n-k_n)},\gamma))^2}\right)\right),$$
so we may represent the ratio of likelihoods in the following form:

\begin{equation}R_n(u) = A_1\cdot A_2\cdot A_3,\label{123}\end{equation} where
$$A_1=\frac{\exp\left(-\sum\limits_{i=0}^{k_n-1}S(X_{(n-i)},\gamma + t(k_n)) +
\sum\limits_{i=0}^{k_n-1}S(X_{(n-i)},\gamma)\right)}{\exp\left(- k_n S(X_{(n-k_n)},\gamma+t(k_n)) + k_n S(X_{(n-k_n)},\gamma) \right)},$$

$$A_2=\left(\frac{S'_x(X_{(n-k_n)},\gamma+t(k_n))}{S'_x(X_{(n-k_n)},\gamma)}\right)^{k_n} \text{ and}$$

$$ A_3=\frac {\left(1 - \frac{S_{xx}^{''}(X_{(n-k_n)},\gamma) }{(S^{'}_x(X_{(n-k_n)},\gamma))^2} + o\left(\frac{S_{xx}^{''}(X_{(n-k_n)},\gamma)}{(S^{'}_x(X_{(n-k_n)},\gamma))^2}\right)\right)^{k_n}}
{\left(1 - \frac{S_{xx}^{''}(X_{(n-k_n)},\gamma+t(k_n)) }{(S^{'}_x(X_{(n-k_n)},\gamma+t(k_n)))^2} + o\left(\frac{S_{xx}^{''}(X_{(n-k_n)},\gamma+t(k_n))}{(S^{'}_x(X_{(n-k_n)},\gamma+t(k_n)))^2}\right)\right)^{k_n}}.$$
Find the asymptotical distribution of $A_1\cdot A_2$. Consider the random variables $\{Y_i\}_{i=1}^{k_n}$ with the $n$th order statistics $Y_{(k_n-i)} = [S(X_{(n-i)},\gamma + t(k_n)) -
S(X_{(n-i)},\gamma)]- [S(X_{(n-k_n)},\gamma+t(k_n)) - S(X_{(n-k_n)},\gamma)], \; i=0, \ldots,
k_n-1,$ that appear in the expression of $A_1.$ According to R\'{e}nyi's representation (for example, see de Haan (2006)), these random variables are independent given $X_{(n-k_n)}=q$. They are also identically distributed. Now consider the multiplier $A_2.$ It is easy to see, that
\[A_2 = \exp( k_n \ln S^{'}_x(q,\gamma+t(k_n)) - k_n \ln S^{'}_x(q,\gamma) )\]
given $X_{(n-k_n)}=q$.
Denote $Z_i = Y_i - \left(\ln S^{'}_x(q,\gamma+t(k_n)) - \ln S^{'}_x(q,\gamma)\right),\ i=1,\ldots, k_n.$ See, that $\sum\limits_{i=1}^{k_n} Z_i = \sum\limits_{i=1}^{k_n} Y_i - \left(\ln S^{'}_x(q,\gamma+t(k_n)) - \ln S^{'}_x(q,\gamma)\right)$ is equal to $-\ln (A_1\cdot A_2)$ in distribution. And it is evident, that $\{Z_i\}_{i=1}^{k_n}$ are independent and identically distributed given $X_{(n-k_n)}=q$ as $\{Y_i\}_{i=1}^{k_n}.$
So, using central limit theorem for triangular arrays we obtain

\begin{equation}\label{N01}
 \frac{\sum\limits_{j=1}^{k_n}Z_j - k_n EZ_1}{\sqrt{k_n DZ_1}} \xrightarrow [k_n\rightarrow\infty]
{d} N(0, 1)\end{equation}
given $X_{(n-k_n)}=q$ and under the assumption that $\lim\limits_{n\to\infty} k_n DZ_1=C_2>0,$ the Lindeberg condition must hold also:
\[M_2(\tau) = k_n E (Y^2_1, |Y_1|>\tau)\to 0\ \ \ as\ n\to\infty\]
for all $\tau>0.$ But it is convenient to verify the Lyapunov condition instead of Lindeberg condition:
\begin{equation}\frac{1}{k_n (DZ_1)^{2}}E(Z_1 - EZ_1)^4\to 0\label{Lya}\end{equation}
as $n\to\infty.$\\
\\
Find the asymptotics of $EZ_1$ и $DZ_1$ given $X_{(n-k_n)}=q.$ Decompose the expressions $S(X_{(n-i)},\gamma + t(k_n)) -
S(X_{(n-i)},\gamma)$ in a Taylor series with the remainder term in Lagrange form for all $i=0,\ldots k_n-1$:

$$S(X_{(n-i)},\gamma + t(k_n)) -
S(X_{(n-i)},\gamma) = t(k_n) S^{'}_{\gamma}(X_{(n-i)},\gamma) +$$
\begin{equation} + \frac{1}{2}(t(k_n))^2 S^{''}_{\gamma \gamma}(X_{(n-i)},\gamma) + \frac{1}{6}(t(k_n))^3 S^{'''}_{\gamma \gamma \gamma}(X_{(n-i)},\gamma + \widetilde{t}(k_n, X_{(n-i)})), \label{teylor}\end{equation}
 where $|\widetilde{t}(k_n, X_{(n-i)})|\leq |t(k_n)|$ and signs of $t(k_n)$ and $\widetilde{t}(k_n, X_{(n-i)})$ are the same. Decompose the expression $\ln S^{'}_x(q,\gamma+t(k_n)) - \ln S^{'}_x(q,\gamma)$ in the same way:

$$\ln S^{'}_x(q,\gamma+t(k_n)) - \ln S^{'}_x(q,\gamma) = t(k_n) \frac{S^{''}_{x\gamma}(q,\gamma)}{S^{'}_x(q,\gamma)} +$$

 \begin{equation}+ \frac{1}{2}(t(k_n))^2\left(\frac{S^{'''}_{x\gamma\gamma}(q,\gamma)}{S^{'}_x(q,\gamma)} - \left(\frac{S^{''}_{x\gamma}(q,\gamma)}{S^{'}_x(q,\gamma)}\right)^2\right) + \frac{1}{6}(t(k_n))^3\left.\left(\frac{S^{''}_{x\gamma}(q,\gamma)}{S^{'}_x(q,\gamma)}\right)^{''}_{\gamma\gamma}\right|_{q, \gamma+\overline{t}(k_n)},\label{teylor2}\end{equation} where $|\overline{t}(k_n)|\leq |t(k_n)|$ and signs of $\overline{t}(k_n)$ and $t(k_n)$ are the same as before.
So, the conditional mathematical expectation of $Z_1$ given $X_{(n-k_n)}=q$ may be represented as following

\[EZ_1 = \frac{\int\limits_{q}^{\infty} [S(x,\gamma + t(k_n)) -
S(x,\gamma)] \exp(-S(x,\gamma)) dx}{ \int\limits_{q}^{\infty} \exp(-S(x,\gamma))
dx} -\]

\[ -[S(q,\gamma + t(k_n)) -
S(q,\gamma)] - [\ln S^{'}_x(q,\gamma+t(k_n)) - \ln S^{'}_x(q,\gamma)] = \]

 \[=t(k_n)\left( \frac{\int\limits_{q}^{\infty}S^{'}_{\gamma}(x,\gamma) \exp(-S(x,\gamma)) dx}{ \int\limits_{q}^{\infty} \exp(-S(x,\gamma))dx} - S^{'}_{\gamma}(q,\gamma) - \frac{S^{''}_{x\gamma}(q,\gamma)}{S^{'}_x(q,\gamma)}\right) + \]

\[+\frac{1}{2}(t(k_n))^2\left( \frac{\int\limits_{q}^{\infty}S^{''}_{\gamma \gamma}(x,\gamma) \exp(-S(x,\gamma)) dx}{ \int\limits_{q}^{\infty} \exp(-S(x,\gamma))dx} - S^{''}_{\gamma \gamma}(q,\gamma) - \frac{S^{'''}_{x\gamma\gamma}(q,\gamma)}{S^{'}_x(q,\gamma)} + \right. \]

\[+ \left.\left(\frac{S^{''}_{x\gamma}(q,\gamma)}{S^{'}_x(q,\gamma)}\right)^2\right) +
\frac{1}{6}(t(k_n))^3\left( \frac{\int\limits_{q}^{\infty}S^{'''}_{\gamma \gamma \gamma}(x,\gamma + \widetilde{t}(k_n, x)) \exp(-S(x,\gamma)) dx}{ \int\limits_{q}^{\infty} \exp(-S(x,\gamma))dx} -\right.\]
\begin{equation} \left.- S^{'''}_{\gamma \gamma \gamma}(q,\gamma + \widehat{t}(k_n)) - \left.\left(\frac{S^{''}_{x\gamma}(x,\gamma)}{S^{'}_x(x,\gamma)}\right)^{''}_{\gamma\gamma}\right|_{q, \gamma+\overline{t}(k_n)}\right), \label{Ey1}\end{equation}
 where $\widehat{t}(k_n) = \widetilde{t}(k_n, q)$. In the sequel the function $S(x,\theta)$ and all of its derivatives will be considered only as $x=q$ and $\theta = \gamma,$ so, for example, the second mixed derivative $S^{''}_{x \gamma}(q,\gamma)$ will be denoted as $S_{x \gamma},$ also we will write $S$ instead of $S(q,\gamma).$ According to Lemma 3.2.1 in de Haan (2006),
$X_{(n-k_n)}=q\to \infty$ as $n\to\infty.$ Find the asymptotics of $I_1 =
\int\limits_{q}^{\infty} \exp(-S) dx.$ Using Lemma 3 we obtain

\begin{equation}I_1 = e^{-S}\left(\frac{1}{S_x} - \frac{S_{xx}}{S_x^3} + \frac{3 S_{xx}^2}{S_x^5} - \frac{S_{xxx}}{S_x^4} + o \left(\max\left(\frac{S_{xx}^2}{S_x^5}, \frac{S_{xxx}}{S_x^4}\right)\right) \right).\label{i1}\end{equation}
 Let's estimate the integral $I_2 = \int\limits_{q}^{\infty} S_{\gamma} \exp(-S) dx,$ it converges according to the regularity condition A3. Using Lemma 2 for the function $S - \ln S_{\gamma}$ (see, that this function satisfies the regularity conditions C1--C3), we obtain

$$I_2 = \frac{S_{\gamma} e^{-S}}{S_x - \frac{S_{x\gamma}}{S_{\gamma}}}\left(1 - \frac{S_{xx} - \frac{S_{xx\gamma}}{S_{\gamma}} + \frac{S^2_{x\gamma}}{S^2_{\gamma}}}{\left(S_x - \frac{S_{x\gamma}}{S_{\gamma}}\right)^2} + \frac{3 \left( S_{xx} - \frac{S_{xx\gamma}}{S_{\gamma}} + \frac{S^2_{x\gamma}}{S^2_{\gamma}}\right)^2}{\left(S_x - \frac{S_{x\gamma}}{S_{\gamma}}\right)^4} \right.- $$

\begin{equation}- \left.\frac{S_{xxx} - \frac{S_{xxx\gamma}}{S_\gamma} + \frac{3S_{x\gamma}S_{xx\gamma}}{S_\gamma^2} - \frac{2S_{x\gamma}^3}{S_\gamma^3}}{\left(S_x - \frac{S_{x\gamma}}{S_{\gamma}}\right)^3} + o \left(\max\left(\frac{S_{xx}^2}{S_x^4}, \frac{S_{xxx}}{S_x^3}\right)\right) \right).\label{i2}\end{equation}
According to the regularity conditions A1--A3 and the L'Hospital rule, we have
\begin{equation}\lim\limits_{q\to\infty}\frac{S_{x\gamma}}{S_x S_{\gamma}} = \lim\limits_{q\to\infty}\frac{S_{x\gamma}/S_\gamma}{S_x} = \lim\limits_{q\to\infty} \frac{\ln S_\gamma}{S} = 0,\label{frac}\end{equation}
consequently,
\begin{equation}\frac{1}{1 - \frac{S_{x\gamma}}{S_x S_{\gamma}}} = 1 + \frac{S_{x\gamma}}{S_x S_{\gamma}} + \frac{S^2_{x\gamma}}{S^2_x S^2_{\gamma}} + o\left(\frac{S^2_{x\gamma}}{S^2_x S^2_{\gamma}} \right).\label{sxg}\end{equation} Here and in the sequel we omit the module brackets in expressions like $\ln S_\gamma,$ since if, for example, $S_\gamma<0,$ then $(\ln (-S_\gamma))_x = \frac{-S_{x\gamma}}{-S_\gamma} = \frac{S_{x\gamma}}{S_\gamma},$ and the limit in (\ref{frac}) does not change.
 Now prove that $(\ln S_\gamma)_{xx} =\frac{S_{xx\gamma}}{S_{\gamma}} - \frac{S^2_{x\gamma}}{S^2_{\gamma}} = o(S_{xx})$ and $(\ln S_\gamma)_{xxx} = \frac{S_{xxx\gamma}}{S_\gamma} - \frac{3S_{x\gamma}S_{xx\gamma}}{S_\gamma^2} + \frac{2S_{x\gamma}^3}{S_\gamma^3} = o(S_{xxx})$ and $q\to\infty.$ According to regularity conditions A1--A3 and the 
 L'Hospital rule, we have

$$\lim\limits_{q\to\infty} \frac{\ln S_\gamma}{\ln S} = \lim\limits_{q\to\infty}\frac{S_{x\gamma}/S_\gamma}{S_x/S} = \lim\limits_{q\to\infty}\frac{S_{x\gamma}S}{S_x S_\gamma} = const,$$
and
\begin{equation}\lim\limits_{q\to\infty} \frac{\ln S_x}{\ln S} = \lim\limits_{q\to\infty}\frac{S_{xx}/S_x}{S_x/S} = \lim\limits_{q\to\infty}\frac{S_{xx}S}{S^2_x} = const.\label{frac3}\end{equation}
Combining two last results we derive

$$\lim\limits_{q\to\infty} \frac{S^2_{x\gamma}S^2}{S^2_{\gamma} S^2_x} \frac{S^2_x}{S_{xx}S} = \lim\limits_{q\to\infty}\frac{S^2_{x\gamma}}{S^2_{\gamma}} \frac{S}{S_{xx}} = const,$$
 i.e. $\frac{S^2_{x\gamma}}{S^2_{\gamma}} = o (S_{xx}).$ According to the regularity conditions A1 and A3, $\lim\limits_{x\to\infty} \frac{\ln S_{x\gamma}}{\ln S} = const>0,$ thus, $\lim\limits_{x\to\infty}\ln S_{x\gamma}=+\infty.$ Using (\ref{frac3}), we derive that $\frac{\ln S_{x\gamma}}{\ln S_\gamma}$ tends to some positive constant as $q\to\infty.$ 
 Using the L'Hospital rule, obtain

$$\lim\limits_{q\to\infty}\frac{S_{xx\gamma}}{S_{\gamma}} \frac{S^2_{\gamma}}{S^2_{x\gamma}} = \lim\limits_{q\to\infty}\frac{S_{xx\gamma}S_{\gamma}}{S^2_{x\gamma}} = \lim\limits_{q\to\infty} \frac{\ln S_\gamma}{\ln S_{x\gamma}} = const,$$
 i.e. $\frac{S_{xx\gamma}}{S_{\gamma}} = O\left(\frac{S^2_{\gamma}}{S^2_{x\gamma}}\right)$ and $(\ln S_\gamma)_{xx} =\frac{S_{xx\gamma}}{S_{\gamma}} - \frac{S^2_{x\gamma}}{S^2_{\gamma}} = o(S_{xx})$ as $q\to\infty.$ That proof of such fact that $(\ln S_\gamma)_{xxx} = \frac{S_{xxx\gamma}}{S_\gamma} - \frac{3S_{x\gamma}S_{xx\gamma}}{S_\gamma^2} + \frac{2S_{x\gamma}^3}{S_\gamma^3} = o(S_{xxx})$ as $q\to\infty$ is analogous. Note as well that $\frac{S_{xx}}{S^2_x} = O\left(\frac{S_{x\gamma}}{S_x S_\gamma}\right).$ Really,

\begin{equation}\lim\limits_{q\to\infty} \frac{S_{xx}}{S^2_x}\frac{S_x S_\gamma}{S_{x\gamma}}  = \lim\limits_{q\to\infty}\frac{S_{xx}/S_x}{S_{x\gamma}/S_\gamma} = \lim\limits_{q\to\infty}\frac{\ln S_x}{\ln S_\gamma}= const.\label{frac2}\end{equation}
 In addition holds $\frac{S_{xxx}}{S^3_x} = O\left(\frac{S^2_{xx}}{S^4_x}\right) = O\left(\frac{S^2_{x\gamma}}{S^2_x S^2_\gamma}\right).$ Similarly to previous reasoning
$$\lim\limits_{q\to\infty}\frac{S_{xxx}}{S^3_x} \frac{S^4_x}{S^2_{xx}} = \lim\limits_{q\to\infty}\frac{S_{xxx}/S_{xx}}{S_{xx}/S_x} = \lim\limits_{q\to\infty} \frac{\ln S_{xx}}{\ln S_{x}} = const.$$
Here such as in other analogous situations the case $\lim\limits_{q\to\infty} \ln S_{xx} = const$ is investigated simply.
Rewrite the integral $I_2$ using obtained results

$$I_2 = \frac{S_{\gamma} e^{-S}}{S_x} \left[1 + \frac{S_{x\gamma}}{S_{\gamma}S_x} + \left(\frac{S_{x\gamma}}{S_{\gamma}S_x}\right)^2 - \frac{S_{xx} - \frac{S_{xx\gamma}}{S_{\gamma}} + \frac{S^2_{x\gamma}}{S^2_{\gamma}}}{S_x^2} - \frac{3S_{xx}S_{x\gamma}}{S_x^3S_\gamma} + \right.$$

$$\left.+ \frac{3 S_{xx}^2 - S_{xxx}S_x}{S_x^4} + o \left(\frac{S^2_{x\gamma}}{S^2_x S^2_\gamma}\right)\right].$$
Find the asymptotics of the ratio $I_2/I_1.$ It is easy to see that

\begin{equation}\frac{I_2}{I_1} = S_{\gamma} \left[1 + \frac{S_{x\gamma}}{S_{\gamma}S_x} + \frac{S_{xx\gamma}}{S_\gamma S_x^2} - \frac{3S_{xx}S_{x\gamma}}{S_x^3S_\gamma} + o \left(\frac{S^2_{x\gamma}}{S^2_x S^2_\gamma}\right)\right].\label{i2i1}\end{equation}
Thus, the first summand in (\ref{Ey1}) is equal to

$$B_1 = t(k_n)\left(\frac{I_2}{I_1} - S_\gamma - \frac{S_{x\gamma}}{S_x}\right) = t(k_n)\left( \frac{S_{xx\gamma}}{S_x^2} - \frac{3S_{xx}S_{x\gamma}}{S_x^3} + o \left(\frac{S^2_{x\gamma}}{S^2_x S_\gamma}\right)\right).$$
Consider the second summand in (\ref{Ey1}). According to the proof of Lemma 2 (see Rodionov (2014)), it holds $F(q)=
\exp(-S(q))\left(c_0 + c_1 + o(c_1)\right)$ for the integral $F(q) = \int\limits_{q}^{+\infty} \exp(-S(x))dx$ as $q\rightarrow +\infty,$ where $c_0 = \left.\frac{1}{S'(x)}\right|_{x=q},$ and $c_1 =
\frac{1}{S'(x)}\frac{d}{dx}\left.\left(\frac{1}{S'(x)}\right)\right|_{x=q}.$ Derive similarly to (\ref{frac}) and (\ref{frac2}), что $\lim\limits_{q\to\infty}\frac{S_{x\gamma\gamma}}{S_x S_{\gamma\gamma}} = 0$ и что $\frac{S_{xx}}{S^2_x} = O\left(\frac{S_{x\gamma\gamma}}{S_x S_{\gamma\gamma}}\right)$ при $q\to\infty.$ Note also that $\frac{S_{xx\gamma\gamma}}{S_{\gamma\gamma}S_x^2} = o\left(\frac{S_{x\gamma\gamma}}{S_{\gamma\gamma}S_x}\right)$ as $q\to\infty,$ since

$$\lim\limits_{q\to\infty} \frac{\frac{S_{xx\gamma\gamma}}{S_{\gamma\gamma}S_x^2}}{\frac{S_{x\gamma\gamma}}{S_{\gamma\gamma}S_x}} = \lim\limits_{q\to\infty} \frac{S_{xx\gamma\gamma}}{S_{x\gamma\gamma}S_x} = \lim\limits_{q\to\infty}\frac{S_{xx\gamma\gamma}/S_{x\gamma\gamma}}{S_x} = \lim\limits_{q\to\infty} \frac{\ln S_{x\gamma\gamma}}{S} = 0.$$
So we can write the integral $I_3 = \int\limits_{q}^{\infty} S_{\gamma\gamma} \exp(-S) dx$ in the following form
$$I_3 = \frac{S_{\gamma\gamma} e^{-S}}{S_x} \left[1 + \frac{S_{x\gamma\gamma}}{S_{\gamma\gamma}S_x} - \frac{S_{xx}}{S^2_x} + O\left(\frac{S^2_{x\gamma\gamma}}{S^2_x S^2_{\gamma\gamma}}\right)\right] \ \ as\ q\to\infty.$$
Hence we obtain the explicit form of the second summand in (\ref{Ey1})
$$B_2 = \frac{1}{2}t^2(k_n)\left(\frac{I_3}{I_1} - S_{\gamma\gamma} - \frac{S_{x\gamma\gamma}}{S_x} + \frac{S^{2}_{x\gamma}}{S^{2}_x}\right) = \frac{1}{2}t^2(k_n)\left(\frac{S^2_{x\gamma}}{S^2_x} + O \left(\frac{S^2_{x\gamma\gamma}}{S^2_x S_{\gamma\gamma}}\right)\right).$$
Using the estimation of the second summand $B_2,$ we may write the third summand in (\ref{Ey1}) in the following form (it should be recalled that $\widetilde{t}(k_n, q) = \widehat{t}(k_n))$
$$B_3 = \frac{1}{6}(t(k_n))^3\left( \frac{\int\limits_{q}^{\infty}S_{\gamma \gamma \gamma}(x,\gamma + \widetilde{t}(k_n, x)) \exp(-S) dx}{ \int\limits_{q}^{\infty} \exp(-S)dx} - S_{\gamma \gamma \gamma}(q,\gamma + \widehat{t}(k_n)) - \right.$$

$$- \left.\left.\left(\frac{S^{''}_{x\gamma}(x,\gamma)}{S^{'}_x(x,\gamma)}\right)_{\gamma\gamma}\right|_{q, \gamma+\overline{t}(k_n)}\right) = \frac{1}{6}(t(k_n))^3\left(\left.\left(\frac{S_{x\gamma\gamma\gamma}}{S_x}\right)\right|_{q, \gamma+\widehat{t}(k_n)}-\right.$$

\begin{equation} - \left.\left.\left(\frac{S_{x\gamma\gamma\gamma}}{S_x} - 3\frac{S_{x\gamma\gamma}S_{x\gamma}}{S_x^2} + 2 \frac{S^3_{x\gamma}}{S^3_x}\right)\right|_{q, \gamma+\overline{t}(k_n)}\right) (1+ o(1)).\label{b3}\end{equation}
 Let $$t(k_n) = \frac{u}{\sqrt{k_n}}\frac{S_x}{S_{x\gamma}},$$ where $u$ is some constant. According to the continuous mapping theorem (see Billingsley (1999)) and Lemma 3, we get
$$\frac{S_x(q, \gamma)}{S_{x\gamma}(q, \gamma)} - \frac{S_x(a_{n/k_n}, \gamma)}{S_{x\gamma}(a_{n/k_n}, \gamma)}\xrightarrow{P}0.$$
Return to the first summand in (\ref{Ey1}) and consider its asymptotics given $t(k_n)$ as $q\to \infty$

$$B_1 = t(k_n)\left( \frac{S_{xx\gamma}}{S_x^2} - \frac{3S_{xx}S_{x\gamma}}{S_x^3}\right)(1+o(1)) = \frac{u}{\sqrt{k_n}}\left( \frac{S_{xx\gamma}}{S_x S_{x\gamma}} - \frac{3S_{xx}}{S_x^2}\right)(1+o(1)).$$
 It follows from (\ref{frac3}), that $\frac{S_{xx}}{S_x^2} = O(S^{-1})$ as $q\rightarrow\infty.$ Prove that $\frac{S_{xx\gamma}}{S_x S_{x\gamma}} = O(S^{-1})$ as $q\rightarrow\infty.$ According to the
  L'Hospital rule and regularity conditions, we have
\begin{equation}\lim\limits_{q\to\infty}\frac{S_{xx\gamma}S}{S_x S_{x\gamma}} = \lim\limits_{q\to\infty} \frac{S_{xx\gamma}/S_{x\gamma}}{S_x/S} = \lim\limits_{q\to\infty}\frac{\ln S_{x\gamma}}{\ln S} = const, \label{frac4}\end{equation}
so $\frac{S_{xx\gamma}}{S_x S_{x\gamma}} = O(S^{-1})$ holds.
  Hence $k_nB_1 = \sqrt{k_n}O(S^{-1})$ as $q\to\infty.$ According to the continuous mapping theorem and Lemma 3, $S(X_{(n-k_n)}, \gamma) - S(a_{n/k_n}, \gamma)\xrightarrow{P} 0,$ while according to the regularity condition A3, $\lim\limits_{a_{n/k_n}\to\infty} \frac{\ln S_x(a_{n/k_n}, \gamma)}{S(a_{n/k_n}, \gamma)} = 0.$ It is obtained in the proof of Lemma 3, that $\frac{\exp(-S(a_{n/k_n},\gamma))}{S_x(a_{n/k_n},\gamma)} =
\frac{k_n}{n}.$ Find the logarithm of both parts of this equation and derive that \begin{equation}\frac{S(X_{(n-k_n)},\gamma)}{\ln \frac{n}{k_n}} \xrightarrow{P} 1.\label{Slon}\end{equation} Since under the condition (\ref{cond}) $\sqrt{k_n} =
o(\ln\frac{n}{k_n}),$ then $k_nB_1 \xrightarrow{P} 0.$ Then it is easy to see using the value of $t(k_n)$ that $k_n B_2 \xrightarrow{P} \frac{u^2}{2}.$ Return to the third summand in (\ref{Ey1}) and prove that $k_nB_3 \xrightarrow{P} 0.$ Firstly note that $\frac{S_{x\gamma\gamma\gamma}}{S_x} = R_1(S),$ $\frac{S_{x\gamma\gamma}S_{x\gamma}}{S_x^2} = R_2(S)$ and $\frac{S_{x\gamma}}{S_x} = R_3(S),$ where $R_1(S),$ $R_2(S)$ and $R_3(S)$ are some slowly varying functions (see for example Galambos and Seneta(1973)). 
From the L'Hospital rule and regularity conditions imply that for all $\gamma>0$

 \[\lim\limits_{q\to\infty}\frac{S_{x\gamma\gamma\gamma}/S_{\gamma\gamma\gamma}}{S_x/S} = \lim\limits_{q\to\infty}\frac{\ln S_{\gamma\gamma\gamma}}{\ln S} = 1,\]

\[\lim\limits_{q\to\infty}\frac{\frac{S_{x\gamma\gamma}S_{x\gamma}}{S_{\gamma\gamma}S_{\gamma}}}{S_x^2/S^2} = \lim\limits_{q\to\infty}\frac{\ln S_{\gamma\gamma} \ln S_\gamma}{(\ln S)^2} = 1,\]

\[\lim\limits_{q\to\infty}\frac{S_{x\gamma}/S_{\gamma}}{S_x/S} = \lim\limits_{q\to\infty}\frac{\ln S_{\gamma}}{\ln S} = 1,\]
 it means that $\frac{S_{x\gamma\gamma\gamma}}{S_x} = O\left(\frac{S_{\gamma\gamma\gamma}}{S}\right),$ $\frac{S_{x\gamma\gamma}S_{x\gamma}}{S_x^2} = O\left(\frac{S_{\gamma\gamma}S_{\gamma}}{S^2}\right)$ and $\frac{S_{x\gamma}}{S_x} = O\left(\frac{S_{\gamma}}{S}\right)$ respectively. But it appears from the regularity condition A3, that $\frac{S_{\gamma\gamma\gamma}}{S},$ $\frac{S_{\gamma\gamma}S_{\gamma}}{S^2}$ are $\frac{S_{\gamma}}{S}$ are slowly varying functions of $S,$ so $\frac{S_{x\gamma\gamma\gamma}}{S_x} = R_1(S),$ $\frac{S_{x\gamma\gamma}S_{x\gamma}}{S_x^2} = R_2(S)$ and $\frac{S_{x\gamma}}{S_x} = R_3(S)$ hold.\\
 \\
 According to Lagrange theorem, $\forall\theta,$ $0<\theta\leq t(k_n),$ there exist such $\widetilde{\theta},$ $0<\widetilde{\theta}\leq\theta$ and the slowly varying function $R(x),$ that \[\left|\frac{S(q,\gamma)}{S(q, \gamma+\theta)}\right| = \left|\frac{S(q,\gamma+\theta) - \theta S_\gamma(q, \gamma + \widetilde{\theta})}{S(q,\gamma+\theta)}\right|\leq \]

$$\leq1 + \theta \left|\frac{S_\gamma(q, \gamma + \widetilde{\theta})}{S(q, \gamma + \widetilde{\theta})}\right| \leq 1 + \frac{1}{\sqrt{k_n}} R(S),$$ since $S$ is strictly monotone starting from some $x_0>0.$ Since under the condition (\ref{cond}), $\lim\limits_{n\to\infty}\frac{k_n}
{(\ln\frac{n}{k_n})^{\varepsilon}} = 1$ for some $\varepsilon,$ $0<\varepsilon<2,$ and, as it is mentioned before, $\frac{S(X_{(n-k_n)},\gamma)}{\ln \frac{n}{k_n}} \xrightarrow{P} 1,$ then $\frac{S(X_{(n-k_n)},\gamma)}{S(X_{(n-k_n)}, \gamma+\theta)}\xrightarrow{P} 1,$ so it holds for the arbitrary slowly varying function $R(x)$ and the arbitrary sequence $\theta(k_n)$ such that $0\leq\theta(k_n)\leq t(k_n),$ that $\frac{R(S(q, \gamma+\theta(k_n)))} {R(S(q, \gamma))}\xrightarrow{P} 1.$ Hence $k_n B_3 = \frac{1}{\sqrt{k_n}} R(S(q, \gamma)) \xrightarrow{P} 0,$ q.e.d.. So, we obtain \begin{equation}k_nEZ_1\xrightarrow{P} \frac{u^2}{2}.\label{E}\end{equation}
 Consider the conditional variance $DZ_1 = EZ_1^2 - (EZ_1)^2$ given $X_{(n-k_n)}=q.$ Note, since $k_nEY_1\xrightarrow{P} u^2/2$ as $n\to\infty,$ then $k_n(EY_1)^2\xrightarrow{P} 0$ as $n\to \infty.$ Find $EZ^2_1$ given $X_{(n-k_n)}=q$

$$ EY^2_1 = \frac{\int\limits_{q}^{\infty} ([S(x,\gamma + t(k_n)) -
S(x,\gamma)] -[S(q,\gamma + t(k_n)) -
S(q,\gamma)] - }{ \int\limits_{q}^{\infty} \exp(-S(x,\gamma))
dx}$$
$$\frac{ - [\ln S^{'}_x(q,\gamma+t(k_n)) - \ln S^{'}_x(q,\gamma)])^2 \exp(-S(x,\gamma)) dx}{}=  $$

$$ = \frac{\int\limits_{q}^{\infty}[S(x,\gamma + t(k_n)) -
S(x,\gamma)]^2 \exp(-S(x,\gamma)) dx}{ \int\limits_{q}^{\infty} \exp(-S(x,\gamma))
dx} - $$
$$ - 2  \frac{\int\limits_{q}^{\infty}[S(x,\gamma + t(k_n)) -
S(x,\gamma)] \exp(-S(x,\gamma)) dx}{ \int\limits_{q}^{\infty} \exp(-S(x,\gamma))
dx}\cdot$$

$$\cdot\left([S(q,\gamma + t(k_n)) -
S(q,\gamma)] - [\ln S^{'}_x(q,\gamma+t(k_n)) - \ln S^{'}_x(q,\gamma)]\right) + $$

\begin{equation} + \left([S(q,\gamma + t(k_n)) -
S(q,\gamma)] - [\ln S^{'}_x(q,\gamma+t(k_n)) - \ln S^{'}_x(q,\gamma)]\right)^2.\label{Dy1}\end{equation}
\normalsize
 Using (\ref{teylor}), represent this expression in polynomial form with respect to $t(k_n)$. The coefficients at the null and the first powers of $t(k_n)$ are equal to zero evidently. Consider the coefficient at $(t(k_n))^2,$ whose asymptotics coincides with the asymptotics of the concerned expression. Similarly to (\ref{i2i1}),

\[\frac{\int\limits_{q}^{\infty}S^{2}_{\gamma} \exp(-S) dx}{ \int\limits_{q}^{\infty} \exp(-S)dx} = S^2_{\gamma} \left[1 + 2\frac{S_{x\gamma}}{S_{\gamma}S_x} + 2\frac{S_{xx\gamma}}{S_\gamma S_x^2} + 2\frac{S^2_{x\gamma}}{S^2_{\gamma}S^2_x} - \frac{6S_{xx}S_{x\gamma}}{S_x^3S_\gamma} + o \left(\frac{S^2_{x\gamma}}{S^2_x S^2_\gamma}\right)\right].\]
Using (\ref{frac3}) and (\ref{frac4}), derive the coefficient at $(t(k_n))^2$:

$$a_2 = \frac{\int\limits_{q}^{\infty}S^{2}_{\gamma} \exp(-S) dx}{ \int\limits_{q}^{\infty} \exp(-S)dx} - 2\frac{\int\limits_{q}^{\infty}S_{\gamma} \exp(-S) dx}{ \int\limits_{q}^{\infty} \exp(-S)dx}\left(S_\gamma + \frac{S_{x\gamma}}{S_x}\right) + \left(S_\gamma + \frac{S_{x\gamma}}{S_x}\right)^2 = $$

$$ = \frac{S^2_{x\gamma}}{S^2_x}\left(1 - 2 \frac{S_{xx\gamma}}{S_{x\gamma}S_x} + 6 \frac{S_{xx}}{S^2_x}\right) +  o \left(\frac{S^2_{x\gamma}}{S^2_x}\right) = \frac{S^2_{x\gamma}}{S^2_x}(1+o(1)).$$ Obtain similarly that the coefficient at $(t(k_n))^3$ is equal to \[a_3 = \frac{S_{x\gamma}S_{x\gamma\gamma}}{S^2_x}(1+o(1)).\] But $a_3 = R(S)a_2,$ where $R(S)$ is some slowly varying function of $S.$ 
From the L'Hospital rule and the regularity condition A3 imply that
 \begin{equation}\lim\limits_{q\to\infty}\frac{S_{x\gamma}/S_{\gamma}}{S_{x\gamma\gamma}/S_{\gamma\gamma}} = \lim\limits_{q\to\infty}\frac{\ln S_{\gamma}}{\ln S_{\gamma\gamma}} = 1,\label{giraf}\end{equation}
 hence $\frac{a_3(t(k_n))^3}{a_2 (t(k_n))^2}\to 0$ as $q\to\infty.$ So, $EZ^2_1 = (t(k_n))^2 \frac{S^2_{x\gamma}}{S^2_x}(1+o(1))$ and
\begin{equation}k_nDZ_1\xrightarrow{P} u^2.\label{D}\end{equation}
 Now verify the Lyapunov condition (\ref{Lya}). See, that $E(Z_1 - EZ_1)^4 = EZ_1^4 - 4EZ_1EZ^3_1 + 6EZ^2_1 (EZ_1)^2 - 3(EZ_1)^4.$ Estimation of $EZ_1^3$ and $EZ_1^4$ is done similarly to estimation of $EY_1^2:$  $EZ_1^3 = (t(k_n))^3 \frac{S^3_{x\gamma}}{S^3_x}(1+o(1))$ and $EZ_1^4 = (t(k_n))^4 \frac{S^4_{x\gamma}}{S^4_x}(1+o(1))$ as $q\to\infty.$ Using derived asymptotics of $EZ_1$ and $EZ_1^2,$ obtain \[E(Z_1 - EZ_1)^4 = (t(k_n))^4 \frac{S^4_{x\gamma}}{S^4_x}(1+o(1)).\] Consequently, we have
$$\frac{1}{k_n (DZ_1)^{2}}E(Z_1 - EZ_1)^4 = \frac{ (t(k_n))^4 \frac{S^4_{x\gamma}}{S^4_x}}{k_n (t(k_n))^4 \frac{S^4_{x\gamma}}{S^4_x}}(1+o(1)) = \frac{1}{k_n}(1+o(1)) \xrightarrow[n\to\infty]{} 0,$$
 so the Lyapunov condition (\ref{Lya}) holds under the conditions of Theorem 1, so the conditions of central limit theorem (\ref{N01}) are satisfied. Since according to Lemma 3, $X_{(n-k_n)} - a_{n/k_n} \xrightarrow{P} 0,$ then it follows from (\ref{E}), (\ref{D}) and the continuous mapping theorem,  that
$$A_1 A_2 \xrightarrow{d} \exp \left( - N\left(\frac{u^2}{2}, u^2\right)\right).$$
Consider the third summand in (\ref{123}) $A_3.$ It is easy to see that
$$\ln A_3 = k_n\left(\frac{S_{xx}(q, \gamma+t(k_n))}{S^2_x(q, \gamma+t(k_n))} - \frac{S_{xx}(q, \gamma)}{S^2_x(q, \gamma)}\right)(1+o(1)),$$ where for some $\theta(k_n),$ $0\leq \theta(k_n) \leq t(k_n)$ holds

$$\frac{S_{xx}(q, \gamma+t(k_n))}{S^2_x(q, \gamma+t(k_n))} - \frac{S_{xx}(q, \gamma)}{S^2_x(q, \gamma)} = $$

$$= t(k_n) \left(\frac{S_{xx\gamma}(q, \gamma+\theta(k_n))}{S^2_x(q, \gamma+\theta(k_n))} - 2 \frac{S_{xx}(q, \gamma+\theta(k_n))S_{x\gamma}(q, \gamma+\theta(k_n))}{S^3_x(q, \gamma+\theta(k_n))}\right).$$
 Consider the asymptotic behavior of $\frac{S_{xx\gamma}}{S_x^2}$ and $\frac{S_{xx}S_{x\gamma}}{S_x^3}$ as $q\to \infty.$ 
 According to the L'Hospital rule and the regularity conditions, we have
$$\lim\limits_{q\to\infty} \frac{S_{xx}S_{x\gamma}S^2}{S_x^3S_\gamma} = \lim\limits_{q\to\infty} \frac{\frac{S_{xx}S_{x\gamma}}{S_x S_\gamma}}{S_x^2/ S^2} = \lim\limits_{q\to\infty} \frac{\ln S_x \ln S_\gamma}{(\ln S)^2} = const,$$
 it means that $\frac{S_{xx}S_{x\gamma}}{S_x^3} = O\left(\frac{S_\gamma}{S^2}\right) = \frac{R(S)}{S},$ where $R(x)$ is some slowly varying function. It is established previously that $\lim\limits_{q\to\infty}\frac{\ln S_{x\gamma}}{\ln S}=const>0,$ so
$$\lim\limits_{q\to\infty}\frac{S_{xx\gamma}}{S_x^2}\frac{S_x^3}{S_{xx}S_{x\gamma}} = \lim\limits_{q\to\infty} \frac{S_{xx}/S_x}{S_{xx\gamma}/S_{x\gamma}} = \lim\limits_{q\to\infty}\frac{\ln S_x}{\ln S_{x\gamma}} = const,$$
 i.e. $\frac{S_{xx\gamma}}{S_x^2} = O\left(\frac{S_{xx}S_{x\gamma}}{S_x^3}\right) = \frac{R'(S)}{S},$ where $R'(x)$ is some slowly varying function. With a glance of previous results and the condition (\ref{cond}) it follows that $\ln A_3 = \sqrt{k_n} \frac{\overline{R}(S)}{S}\xrightarrow{P} 0,$ where $\overline{R}(x)$ is some slowly varying function again, and
$$A_3 \xrightarrow{P} 1.$$
Hence, according to Slutsky's theorem,

$$R_n(u)\xrightarrow{d} \exp\left(-N\left(\frac{u^2}{2}, u^2\right)\right).$$
The proof of Theorem 1 is complete. $\blacksquare$\\
\\
{\bf The proof of Theorem 2.}\\
\\
The scheme of the proof of Theorem 2 is the same as the proof of Theorem 1. Using Lemma 1, write the ratio of likelihoods
\[R_n(t) = \frac{\prod\limits_{i=0}^{k_n-1} \exp[-S(X_{(n-i)},\gamma + t(k_n)) +
S(X_{(n-i)},\gamma)]} {\left ( \int\limits_{X_{(n-k_n)}}^{+\infty}
\exp[-S(x,\gamma)] dx\right)^{-k_n}\cdot
\left(\int\limits_{X_{(n-k_n)}}^{+\infty}
\exp[-S(x,\gamma+t(k_n))]dx \right)^{k_n}}.\]
 Since the concerned family of densities satisfies the conditions of Lemma 2, then the expansion (\ref{123}) holds for the ratio of likelihoods given $X_{(n-k_n)} = q$. As before, consider the random variables$\{Y_i\}_{i=1}^{k_n}$ with the following $n$th order statistics $$Y_{(k_n-i)} = [S(X_{(n-i)},\gamma + t(k_n)) -
S(X_{(n-i)},\gamma)]- $$
$$- [S(X_{(n-k_n)},\gamma+t(k_n)) - S(X_{(n-k_n)},\gamma)],$$ where $i=0, \ldots,
k_n-1,$ that are independent given $X_{(n-k_n)}=q$ according to R\'{e}nyi's representation and identically distributed in addition. It is easy to see that $\ln A_1 = - \sum\limits_{i=1}^{k_n} Y_i.$ Then
$$\ln A_2 =  k_n \ln S^{'}_x(q, \gamma+t(k_n)) - k_n \ln S^{'}_x(q, \gamma).$$
It follows from the L'Hospital rule and the regularity condition B4, that it holds for all $\gamma$ that $S(x, \gamma)$ satisfies the regularity conditions B1-B4

$$\lim\limits_{q\to\infty} \frac{S^{''}_{xx}(q, \gamma)}{(S^{'}_x(q, \gamma))^2} = \lim\limits_{q\to\infty} \frac{\frac{S^{''}_{xx}(q, \gamma)}{S^{'}_x(q, \gamma)}}{S^{'}_x(q, \gamma)} = \lim\limits_{q\to\infty} \frac {\ln S^{'}_x(q, \gamma)}{S(q, \gamma)} = 0,$$
so we have
$$\ln A_3 = k_n\left(\frac{S^{''}_{xx}(q, \gamma+t(k_n))}{(S^{'}_x(q, \gamma+t(k_n)))^2} - \frac{S^{''}_{xx}(q, \gamma)}{(S^{'}_x(q, \gamma))^2}\right)(1+o(1)).$$
Denote
$$G(q) = \left(\ln S^{'}_x(q, \gamma+t(k_n)) - \ln S^{'}_x(q, \gamma)\right) + $$
$$+ \left(\frac{S^{''}_{xx}(q, \gamma+t(k_n))}{(S^{'}_x(q, \gamma+t(k_n)))^2} - \frac{S^{''}_{xx}(q, \gamma)}{(S^{'}_x(q, \gamma))^2}\right),$$

$$H(q) = \sqrt{k_n}t(k_n)\left(\frac{\int\limits_{q}^{\infty}S^{'}_{\gamma}(x,\gamma) \exp(-S(x,\gamma)) dx}{ \int\limits_{q}^{\infty} \exp(-S(x,\gamma))dx} - S^{'}_{\gamma}(q,\gamma) - \right.$$

$$\left.- \frac{S^{''}_{x\gamma}(q,\gamma)}{S^{'}_x(q,\gamma)} - \frac{S^{'''}_{xx\gamma}(q, \gamma)}{(S^{'}_x(q, \gamma))^2} + 2 \frac{S^{''}_{xx}(q, \gamma)S^{''}_{x\gamma}(q, \gamma)}{(S^{'}_x(q, \gamma))^3}\right).$$
 Consider the random variables $\{Z_i\}_{i=1}^{k_n},$ where $Z_i = Y_i - G(q) - \frac{1}{\sqrt{k_n}}H(q).$ It is easy to see that these random variables are identically distributed and independent given $X_{(n-k_n)} = q.$ Note also that $$\sum\limits_{i=1}^{k_n} Z_i = - \ln R_n(t) + \sqrt{k_n} H(q).$$
According to Lindeberg central limit theorem, we have

\begin{equation}\label{N012}
 \frac{\sum\limits_{j=1}^{k_n}Z_j - k_n EZ_1}{\sqrt{k_n DZ_1}} \xrightarrow [k_n\rightarrow\infty]
{d} N(0, 1),\end{equation}
on condition that $X_{(n-k_n)}=q,$ $\lim\limits_{n\to\infty} k_n DZ_1=const>0$ and the Lindeberg condition holds. As before instead of Lindeberg condition we verify Lyapunov condition that take on following form:
\begin{equation}\frac{1}{k_n (DZ_1)^{2}}E(Z_1 - EZ_1)^4\to 0\label{Lya2}.\end{equation}
as $k_n\to\infty.$
Find the asymptotics of $EZ_1$ and $DZ_1$ given $X_{(n-k_n)}=q.$ Firstly note that

$$\frac{S^{''}_{xx}(q, \gamma+t(k_n))}{(S^{'}_x(q, \gamma+t(k_n)))^2} - \frac{S^{''}_{xx}(q, \gamma)}{(S^{'}_x(q, \gamma))^2} = t(k_n) \left(\frac{S^{'''}_{xx\gamma}(q, \gamma)}{(S^{'}_x(q, \gamma))^2} - \right.$$

\begin{equation}\left. - 2 \frac{S^{''}_{xx}(q, \gamma)S^{''}_{x\gamma}(q, \gamma)}{(S^{'}_x(q, \gamma))^3}\right) + \frac{(t(k_n))^2}{2} \left.\left(\frac{S^{''}_{xx}(x, \gamma)}{(S^{'}_x(x, \gamma)^2}\right)^{''}_{\gamma\gamma}\right|_{(x, \gamma) = (q, \gamma+\theta(k_n))},\label{teylor3}\end{equation}
 where $|\theta(k_n)|\leq |t(k_n)|$ and the signs of $\theta(k_n)$ and $t(k_n)$ are the same. It appears from (\ref{teylor}), (\ref{teylor2}) and the previous expansion, that:

$$EZ_1 = \frac{\int\limits_{q}^{\infty} [S(x,\gamma + t(k_n)) -
S(x,\gamma)] \exp(-S(x,\gamma)) dx}{ \int\limits_{q}^{\infty} \exp(-S(x,\gamma))
dx} - $$

$$-[S(q,\gamma + t(k_n)) -
S(q,\gamma)] - G(q) - \frac{1}{\sqrt{k_n}}H(q) = $$

$$ =t(k_n)\left( \frac{\int\limits_{q}^{\infty}S^{'}_{\gamma}(x,\gamma) \exp(-S(x,\gamma)) dx}{ \int\limits_{q}^{\infty} \exp(-S(x,\gamma))dx} - S^{'}_{\gamma}(q,\gamma) - \frac{S^{''}_{x\gamma}(q,\gamma)}{S^{'}_x(q,\gamma)} - \right.$$

$$\left.- \frac{S^{'''}_{xx\gamma}(q, \gamma)}{(S^{'}_x(q, \gamma))^2} + 2 \frac{S^{''}_{xx}(q, \gamma)S^{''}_{x\gamma}(q, \gamma)}{(S^{'}_x(q, \gamma))^3} - \frac{H(q)}{\sqrt{k_n}t(k_n)}\right) + $$

$$+ \frac{1}{2}(t(k_n))^2\left( \frac{\int\limits_{q}^{\infty}S^{''}_{\gamma \gamma}(x,\gamma) \exp(-S(x,\gamma)) dx}{ \int\limits_{q}^{\infty} \exp(-S(x,\gamma))dx} - S^{''}_{\gamma \gamma}(q,\gamma) - \frac{S^{'''}_{x\gamma\gamma}(q,\gamma)}{S^{'}_x(q,\gamma)} + \right.$$

$$+ \left. \left(\frac{S^{''}_{x\gamma}(q,\gamma)}{S^{'}_x(q,\gamma)}\right)^2 - \left.\left(\frac{S^{''}_{xx}(x, \gamma)}{S^{'}_x(x, \gamma)^2}\right)^{''}_{\gamma\gamma}\right|_{(q, \gamma+\theta(k_n))} \right) + $$

$$+ \frac{1}{6}(t(k_n))^3\left( \frac{\int\limits_{q}^{\infty}S^{'''}_{\gamma \gamma \gamma}(x,\gamma + \widetilde{t}(k_n, x)) \exp(-S(x,\gamma)) dx}{ \int\limits_{q}^{\infty} \exp(-S(x,\gamma))dx} - \right.$$

\begin{equation}\left. - S^{'''}_{\gamma \gamma \gamma}(q,\gamma + \widehat{t}(k_n)) - \left.\left(\frac{S^{''}_{x\gamma}(x,\gamma)}{S^{'}_x(x,\gamma)}\right)^{''}_{\gamma\gamma}\right|_{q, \gamma+\overline{t}(k_n)}\right), \label{Ey12}\end{equation}
 where $\widehat{t}(k_n) = \widetilde{t}(k_n, q)$. We omit the arguments of the function $S(x, \gamma)$ and all
  of its partial derivatives as before provided that are equal to $q$ and $\gamma$ respectively.
  It appears from the definition of the function $H(x)$ that the term at $t(k_n)$
  in the expansion (\ref{Ey12}) is equal to 0 identically. Consider the term at $(t(k_n))^2$ in the expansion (\ref{Ey12}).
  Recall that it holds for the integral $F(q) = \int\limits_{q}^{\infty} \exp(-S(x))dx$ as $q\to +\infty$ the following
 \[F(q)=
\exp(-S(q))\left(c_0 + c_1 + o(c_1)\right),\] where $c_0 = \left.\frac{1}{S'(x)}\right|_{x=q}$ and $c_1 =
\frac{1}{S'(x)}\frac{d}{dx}\left.\left(\frac{1}{S'(x)}\right)\right|_{x=q}.$ Using this fact, obtain:

$$I_1 =  \int\limits_{q}^{\infty} \exp(-S)dx = \frac{e^{-S}}{S_x}\left(1 - \frac{S_{xx}}{S_x^2} + o\left(\frac{S_{xx}}{S_x^2}\right)\right),$$

\[I_3 = \int\limits_{q}^{\infty} S_{\gamma\gamma}\exp(-S)dx = \]

$$= \frac{e^{-S}}{S_x - \frac{S_{x\gamma\gamma}}{S_{\gamma\gamma}}}
\left(1 - \frac{S_{xx} - \frac{S_{xx\gamma\gamma}}{S_{\gamma\gamma}}
 + \frac{S^2_{x\gamma\gamma}}{S^2_{\gamma\gamma}}}{\left(S_x - \frac{S_{x\gamma\gamma}}
 {S_{\gamma\gamma}}\right)^2} + o\left(\frac{S_{xx} - \frac{S_{xx\gamma\gamma}}{S_{\gamma\gamma}}
  + \frac{S^2_{x\gamma\gamma}}{S^2_{\gamma\gamma}}}{\left(S_x - \frac{S_{x\gamma\gamma}}{S_{\gamma\gamma}}\right)^2}\right)\right).$$
Under the regularity condition B4 we have
$\frac{S_{x\gamma\gamma}}{S_x S_{\gamma\gamma}} = O\left(\frac{S_{xx}}{S_x^2}\right)$ as $q\to +\infty.$
Really, if $S_{xx}$ is not equal to 0 identically in some neighbourhood of infinity,
then it appears from the L'Hospital rule and the regularity conditions B1--B4 that
\[\lim\limits_{q\to\infty} \frac{S_{x\gamma\gamma}}{S_x S_{\gamma\gamma}}
\frac{S^2_x}{S_{xx}}  = \lim\limits_{q\to\infty}\frac{S_{x\gamma\gamma}/
S_{\gamma\gamma}}{S_{xx}/S_x} = \lim\limits_{q\to\infty}\frac{\ln S_{\gamma\gamma}}{\ln S_x} = const.\]
It follows from the proof of Lemma 2 that $\lim\limits_{q\to\infty} \frac{S_{xx}}{S^2_x} = 0,$
so $\lim\limits_{q\to\infty} \frac{S_{x\gamma\gamma}}{S_x S_{\gamma\gamma}} = 0$ and we get

$$\frac{1}{1 - \frac{S_{x\gamma\gamma}}{S_x S_{\gamma\gamma}}} = 1 + \frac{S_{x\gamma\gamma}}{S_x S_{\gamma\gamma}} + o\left(\frac{S_{x\gamma\gamma}}{S_x S_{\gamma\gamma}} \right).$$
Now prove that $\frac{S_{xx\gamma\gamma}}{S_{\gamma\gamma}S^2_x}=
 o\left(\frac{S_{xx}}{S^2_x}\right)$ as $q\to +\infty.$ From the L'Hospital rule and the regularity conditions imply that
\[\lim\limits_{q\to\infty} \frac{\frac{S_{xx\gamma\gamma}}{S_{\gamma\gamma}S_x^2}}
{\frac{S_{x\gamma\gamma}}{S_{\gamma\gamma}S_x}} = \lim\limits_{q\to\infty}
\frac{S_{xx\gamma\gamma}}{S_{x\gamma\gamma}S_x} = \lim\limits_{q\to\infty}
\frac{S_{xx\gamma\gamma}/S_{x\gamma\gamma}}{S_x} = \lim\limits_{q\to\infty} \frac{\ln S_{x\gamma\gamma}}{S} = 0,\]
so $\frac{S_{xx\gamma\gamma}}{S_{\gamma\gamma}S^2_x} =
o\left(\frac{S_{x\gamma\gamma}}{S_x S_{\gamma\gamma}}\right) =
 o\left(\frac{S_{xx}}{S^2_x}\right),$ q.e.d. Thus $\frac{S_{xx\gamma\gamma}}{S_{\gamma\gamma}} -
 \frac{S^2_{x\gamma\gamma}}{S^2_{\gamma\gamma}} = o(S_{xx})$ as $q\to +\infty.$
Using previous facts and Lemma 2 we conclude that
$$\frac{I_3}{I_1} = S_{\gamma\gamma} \left(1 + \frac{S_{x\gamma\gamma}}{S_x S_{\gamma\gamma}} + o\left(\frac{S_{xx}}{S_x^2}\right)\right).$$
So the term at $(t(k_n))^2$ in the expansion (\ref{Ey12}) is equal to
$$\frac{(t(k_n))^2}{2}\left(\frac{I_3}{I_1} - S_{\gamma\gamma} - \frac{S_{x\gamma\gamma}}{S_x} + \frac{S^{2}_{x\gamma}}{S^{2}_x} - \left.\left(\frac{S_{xx}}{S^2_x}\right)_{\gamma\gamma}\right|_{(q, \gamma+\theta(k_n))}\right) = $$

$$ = \frac{(t(k_n))^2}{2}\left(\frac{S^{2}_{x\gamma}}{S^{2}_x} - \left.\left(\frac{S_{xx}}{S^2_x}\right)_{\gamma\gamma}\right|_{(q, \gamma+\theta(k_n))} + o\left(\frac{S_{xx}}{S_x^2}\right)\right). $$
Analyze the third summand in the expansion (\ref{Ey12}) (denote it as $B_3$) similarly to (\ref{b3})
$$B_3 = \frac{1}{6}(t(k_n))^3\left(\left.\left(\frac{S_{x\gamma\gamma\gamma}}{S_x}\right)\right|_{q, \gamma+\widehat{t}(k_n)}- \right.$$

$$\left.- \left.\left(\frac{S_{x\gamma\gamma\gamma}}{S_x} - 3\frac{S_{x\gamma\gamma}S_{x\gamma}}{S_x^2} + 2 \frac{S^3_{x\gamma}}{S^3_x}\right)\right|_{q, \gamma+\overline{t}(k_n)}\right) (1+ o(1)).$$
Let
$$t(k_n) = \frac{u}{\sqrt{k_n}}\frac{S_x}{S_{x\gamma}},$$ where $u$ is some constant. Note, that the state of Lemma 3 holds as before,
i.e. $X_{(n-k_n)} - a_{n/k_n} \xrightarrow{P} 0$ as $n\to +\infty.$ See that from the regularity conditions imply
 $\lim\limits_{x\to+\infty}\frac{\ln S_x(x, \gamma)}{S(x, \gamma)} = 0.$ Since it follows from the proof of Lemma 3 that
  $S(a_{n/k_n}, \gamma) + \ln S_x(a_{n/k_n}, \gamma) = \ln \frac{n}{k_n} (1+o(1)).$
   Then $\forall \; \varepsilon>0$ $$\lim\limits_{x\to+\infty}S_x(a_{n/k_n}, \gamma)\left(\frac{n}{k_n}\right)^{\varepsilon} = +\infty.$$
So using regularity condition B1 and (\ref{cond2}) we obtain that $\lim\limits_{x\to+\infty} \sqrt{k_n} S_x(a_{n/k_n}, \gamma) = +\infty,$
 whence it appears $X_{(n-k_n)} - a_{n/k_n} \xrightarrow{P} 0$.\\
\\
From the continuous mapping theorem imply that
\[\frac{S_x(q, \gamma)}{S_{x\gamma}(q, \gamma)} - \frac{S_x(a_{n/k_n}, \gamma)}{S_{x\gamma}(a_{n/k_n}, \gamma)}\xrightarrow{P}0\] as before.
Find the asymptotics of $k_nEZ_1$ as $q\to \infty$ and given $t(k_n).$ Consider the second summand in (\ref{Ey12})
$$k_nB_2 = \frac{u^2}{2}\frac{S^2_x}{S^2_{x\gamma}}\left(\frac{S^{2}_{x\gamma}}{S^{2}_x}
 - \left.\left(\frac{S_{xx}}{S^2_x}\right)_{\gamma\gamma}\right|_{(q, \gamma+\theta(k_n))} + o\left(\frac{S_{xx}}{S_x^2}\right)\right).$$
It is easy to see that
$$\left(\frac{S_{xx}}{S^2_x}\right)_{\gamma\gamma} = \frac{S_{xx\gamma\gamma}}{S_x^2} - 4\frac{S_{xx\gamma}S_{x\gamma}}{S_x^3}
- 2 \frac{S_{xx}S_{x\gamma\gamma}}{S_x^3} + 6 \frac{S_{xx}S^2_{x\gamma}}{S_x^4}.$$
Note that it holds for all such $\gamma$ that $S(x,\gamma)$ is defined:
$\frac{S_{xx\gamma\gamma}}{S_x^2} = O\left(\frac{S_{xx}S_{x\gamma\gamma}}{S_x^3}\right)$ and
 $\frac{S_{xx\gamma}S_{x\gamma}}{S_x^3} = O\left(\frac{S_{xx}S^2_{x\gamma}}{S_x^4}\right)$ as $q\to\infty.$
 Prove these facts. Using the L'Hospital rule and the regularity conditions, we derive

$$\lim\limits_{q\to\infty}\frac{S_{xx\gamma\gamma}}{S_x^2}\frac{S_x^3}{S_{xx}S_{x\gamma\gamma}}
= \lim\limits_{q\to\infty}\frac{S_{xx\gamma\gamma}/S_{x\gamma\gamma}}{S_{xx}/S_x} =
\lim\limits_{q\to\infty}\frac{\ln S_{x\gamma\gamma}}{\ln S_x} = const,$$

$$\lim\limits_{q\to\infty}\frac{S_{xx\gamma}S_{x\gamma}}{S_x^3}\frac{S_x^4}{S_{xx}S^2_{x\gamma}} = \lim\limits_{q\to\infty}\frac{S_{xx\gamma}/S_{x\gamma}}{S_{xx}/S_x} = \lim\limits_{q\to\infty} \frac{\ln S_{x\gamma}}{\ln S_x} = const,$$
whence it appears the concerned facts. Then we have
$$\lim\limits_{q\to\infty}\frac{S_{xx}S_{x\gamma\gamma}}{S_x^3} \frac{S_x^4}{S_{xx}S^2_{x\gamma}} = \lim\limits_{q\to\infty} \frac{S_{x\gamma\gamma}S_x}{S^2_{x\gamma}} = \lim\limits_{q\to\infty}\frac{S_{x\gamma\gamma}/S_{\gamma\gamma}}{S_{x\gamma}/S_\gamma} \lim\limits_{q\to\infty}\frac{S_{x}/S}{S_{x\gamma}/S_\gamma} \lim\limits_{q\to\infty}\frac{S_{\gamma\gamma}S}{S^2_\gamma} = $$

$$=\lim\limits_{q\to\infty}\frac{\ln S_{\gamma\gamma}}{\ln S_\gamma}
\lim\limits_{q\to\infty}\frac{\ln S}{\ln S_\gamma}\lim\limits_{q\to\infty}\frac{S_{\gamma\gamma}S}{S^2_\gamma} = \lim\limits_{q\to\infty}\frac{S_{\gamma\gamma}S}{S^2_\gamma},$$
hence and from the regularity condition B4 imply that $\frac{S_{xx}S_{x\gamma\gamma}}{S_x^3} = \frac{S_{xx}S^2_{x\gamma}}{S_x^4} R(S),$
where $R(x)$ is some slowly varying function. So, we obtain that $$\left.\left(\frac{S_{xx}}{S^2_x}\right)_{\gamma\gamma}\right|_{(q, \gamma+\theta(k_n))} = \left.\frac{S_{xx}S^2_{x\gamma}}{S_x^4} R(S)\right|_{(q, \gamma+\theta(k_n))}.$$
As it is mentioned previously, $\frac{S_{x\gamma}}{S_x} = R(S)$ and $R(S(q, \gamma)) = R_1(S(q, \gamma+\theta(k_n)))$ because of continuity
 of the function $S(x, \gamma)$ with respect to $\gamma,$ where $R(x)$ and $R_1(x)$ are some slowly varying functions. So we obtain
$$k_n (t(k_n))^2\left.\left(\frac{S_{xx}}{S^2_x}\right)_{\gamma\gamma}\right|_{(q, \gamma+\theta(k_n))} =
 \left.\frac{S_{xx}}{S^2_x} R_2(S)\right|_{(q, \gamma+\theta(k_n))},$$
where $R_2(x)$ is some slowly varying function.
Recall that according to the regularity condition B2 there exists $\delta>0$ such that $\lim\limits_{q\rightarrow +\infty}
\frac{\ln S_x}{S^{1-\delta}} = 0,$ hence imply that $\frac{S_{xx}}{S^2_x} = o(S^{-\delta}),$ since

$$\lim\limits_{q\to\infty}\frac{\ln S_x}{S^{1-\delta}} = \lim\limits_{q\to\infty}\frac{S_{xx}/S_x}{S_x S^{-\delta}}
 = \lim\limits_{q\to\infty}\frac{S_{xx}}{S^2_x} S^{\delta} = 0.$$
So, we obtain that
$$k_n (t(k_n))^2\left.\left(\frac{S_{xx}}{S^2_x}\right)_{\gamma\gamma}\right|_{(q, \gamma+\theta(k_n))} \xrightarrow{P} 0$$
as $k_n\to\infty$ and finally
$$k_nB_2 \xrightarrow{P} \frac{u^2}{2}.$$
The proof of the following fact $k_nB_3 \xrightarrow{P} 0,$ where $B_3$ is the third summand in the expansion (\ref{Ey12}),
agrees completely with the proof of the analogous fact in the previous section. So, we obtain
\begin{equation}k_nEZ_1\xrightarrow{P} \frac{u^2}{2}.\label{E2}\end{equation}
Now consider the conditional variance $DZ_1 = EZ_1^2 - (EZ_1)^2$ given $X_{(n-k_n)} = q.$ Note as before that
since $k_nEZ_1\xrightarrow{P} u^2/2$ as $n\to\infty,$ then $k_n(EZ_1)^2\xrightarrow{P} 0$ as $n\to \infty.$ Denote

$$J(q) = G(q) + \frac{1}{\sqrt{k_n}}H(q) + [S(q, \gamma+t(k_n)) - S(q, \gamma)],$$
so we get $Z_{(k_n-i)} = [S(X_{(n-i)},\gamma + t(k_n)) -
S(X_{(n-i)},\gamma)] - J(X_{(n-k_n)}),$ where $i=0, \ldots, k_n-1.$ Find the asymptotics of $EZ^2_1$ as $q\to\infty$

$$ EZ^2_1 = \frac{\int\limits_{q}^{\infty} \left([S(x,\gamma + t(k_n)) -
S(x,\gamma)] - J\right)^2 \exp(-S(x,\gamma)) dx}{ \int\limits_{q}^{\infty} \exp(-S(x,\gamma))
dx} =  $$

$$ = \frac{\int\limits_{q}^{\infty}[S(x,\gamma + t(k_n)) -
S(x,\gamma)]^2 \exp(-S(x,\gamma)) dx}{ \int\limits_{q}^{\infty} \exp(-S(x,\gamma))
dx} - $$

\begin{equation} - 2\frac{\int\limits_{q}^{\infty}[S(x,\gamma + t(k_n)) -
S(x,\gamma)] \exp(-S(x,\gamma)) dx}{ \int\limits_{q}^{\infty} \exp(-S(x,\gamma))
dx}\cdot J(q) + J^2(q).\label{Dy12}\end{equation}
\normalsize
Firstly note that $$(t(k_n))^3 \frac{\int\limits_{q}^{\infty}S_{\gamma\gamma\gamma}
(x, \gamma+\widetilde{t}(k_n, x)) \exp(-S(x, \gamma)) dx}{\int\limits_{q}^{\infty} \exp(-S(x, \gamma))dx} = $$

$$= (t(k_n))^3 S_{\gamma\gamma\gamma}(q, \gamma+\widehat{t}(k_n))(1+o(1)) = \frac{1}{(k_n)^{3/2}} S(q, \gamma+\widehat{t}(k_n)) R(S),$$
where as before $|\widehat{t}(k_n)|\leq |t(k_n)|$ and the signs of $\widehat{t}(k_n)$ and $t(k_n)$ are the same, $R(x)$ is some
slowly varying function. From (\ref{Slon}) imply
 $\frac{S(X_{(n-k_n)},\gamma)}{\ln \frac{n}{k_n}} \xrightarrow{P} 1.$
Hence it appears from (\ref{cond2}) and continuity of the function $S(x, \gamma)$ in respect to $\gamma$, that
 $\frac{1}{\sqrt{k_n}} S(q, \gamma+\widehat{t}(k_n)) R(S) \xrightarrow{P} 0.$ It means that

$$(t(k_n))^3 \frac{\int\limits_{q}^{\infty}S_{\gamma\gamma\gamma}(x, \gamma+\widetilde{t}(k_n, x))
\exp(-S(x, \gamma)) dx}{\int\limits_{q}^{\infty} \exp(-S(x, \gamma))dx} = o\left(\frac{1}{k_n}\right),$$
as $n\to \infty$ and $q\to\infty.$
Hence and from the expansions (\ref{teylor}), (\ref{teylor2}) and (\ref{teylor3}) imply:

$$ J(q) = t(k_n) \frac{\int\limits_{q}^{\infty}S_{\gamma} \exp(-S) dx}{\int\limits_{q}^{\infty} \exp(-S)dx}
 + \frac{t^2(k_n)}{2} \left(\frac{\int\limits_{q}^{\infty}S_{\gamma\gamma}
  \exp(-S) dx}{\int\limits_{q}^{\infty} \exp(-S)dx} + \frac{S^2_{x\gamma}}{S_x^2}\right) + o\left(\frac{1}{k_n}\right).$$
Using the last expansion and the expansion (\ref{teylor}), represent $EZ_1^2$ in polynomial form in respect to $t(k_n).$
It is easy to see, that the coefficients at the null and the first powers of $t(k_n)$ are equal to 0.
Consider the coefficient at $(t(k_n))^2$ whose asymptotics agrees with asymptotics of $EZ_1^2$.
Recall that under the conditions of Theorem 2 holds $\frac{S_{x\gamma}}{S_x S_\gamma} = O\left(\frac{S_{xx}}{S^2_x}\right),$
 so obtain similarly to (\ref{i2i1})

$$\frac{\int\limits_{q}^{\infty}S^{2}_{\gamma} \exp(-S) dx}{ \int\limits_{q}^{\infty} \exp(-S)dx} = $$

$$ = S^2_{\gamma} \left[1 + 2\frac{S_{x\gamma}}{S_{\gamma}S_x} + 2\frac{S_{xx\gamma}}{S_\gamma S_x^2}
 + 2\frac{S^2_{x\gamma}}{S^2_{\gamma}S^2_x} - \frac{6S_{xx}S_{x\gamma}}{S_x^3S_\gamma} + o \left(\frac{S^2_{xx}}{S^4_x}\right)\right].$$
So the coefficient at $(t(k_n))^2$ (denote it as $a_2)$ is equal to

$$a_2 = \frac{\int\limits_{q}^{\infty}S^{2}_{\gamma} \exp(-S) dx}{ \int\limits_{q}^{\infty} \exp(-S)dx}
 - \left(\frac{\int\limits_{q}^{\infty}S_{\gamma} \exp(-S) dx}{ \int\limits_{q}^{\infty} \exp(-S)dx}\right)^2
  = \frac{S^2_{x\gamma}}{S^2_x} - 2 \frac{S_{x\gamma}S_{xx\gamma}}{S_x^3} - \frac{S^2_{xx\gamma}}{S^4_x}.$$
Show that $\frac{S_{xx\gamma}}{S^2_x} = o\left(\frac{S_{x\gamma}}{S_x}\right)$ as $q\to\infty.$
It appears from the L'Hospital rule and the regularity conditions that

\begin{equation}\lim\limits_{q\to\infty}\frac{S_{xx\gamma}}{S^2_x}\frac{S_x}{S_{x\gamma}}
 = \lim\limits_{q\to\infty}\frac{S_{xx\gamma}/S_{x\gamma}}{S_x} = \lim\limits_{q\to\infty}\frac{\ln S_{x\gamma}}{S}
  = 0,\label{pavian}\end{equation}
hence the concerned fact appears. Thus we obtain that

$$a_2 = \frac{S^2_{x\gamma}}{S^2_x}(1+o(1)).$$
Similarly, the coefficient at $(t(k_n))^3$ is equal to $a_3 = \frac{S_{x\gamma}S_{x\gamma\gamma}}{S^2_x}(1+o(1)).$
But it follows from (\ref{giraf}) that $a_3 = R(S)a_2,$
so we obtain \[\frac{a_3(t(k_n))^3}{a_2 (t(k_n))^2} = \frac{R_1(S)}{\sqrt{k_n}}\xrightarrow{P} 0,\] where $R(S)$ and $R_1(S)$
are some slowly varying functions. Finally we get $EZ^2_1 = (t(k_n))^2 \frac{S^2_{x\gamma}}{S^2_x}(1+o(1))$ and consequently
\begin{equation}k_nDZ_1\xrightarrow{P} u^2.\label{D2}\end{equation}
Verification of Lyapunov condition (\ref{Lya2}) is done similarly to verification of the condition (\ref{Lya}).
Thus the conditions of the central limit theorem (\ref{N012}) holds.
Since it follows from Lemma 3 that $X_{(n-k_n)} - a_{n/k_n} \xrightarrow{P} 0,$ then
from (\ref{E2}), (\ref{D2}) and the continuous mapping theorem imply
$$R_n(u)\exp(-\sqrt{k_n} H(q)) \xrightarrow{d} \exp \left( - N\left(\frac{u^2}{2}, u^2\right)\right).$$
Now consider the asymptotic behavior of $\exp(\sqrt{k_n} H(q))$ as $n\to\infty.$
As it is stated in Lemma 3,
$$\sqrt{k_n}S_x(a_{n/k_n},\gamma)(X_{(n-k_n)} - a_{n/k_n}) \xrightarrow {d} N(0,1).$$
Denote $\widetilde{H}(x) = \frac{H(x)}{S_x(x,\gamma)},$ $b_n = \frac{1}{\sqrt{k_n}S_x(a_{n/k_n},\gamma)},$ and
 $\eta_n = \sqrt{k_n}S_x(a_{n/k_n},\gamma)(X_{(n-k_n)} - a_{n/k_n}).$
Then according to continuity theorem (see, for example, theorem 1.5.3 in Borovkov(1984)),
if there exists the finite limit of the derivative of the function $\widetilde{H}(x)$ as $x\to\infty,$ then

$$ \frac{\widetilde{H}(a_{n/k_n} + b_n \eta_n) - \widetilde{H}(a_{n/k_n})}{b_n}
\xrightarrow{d} \lim\limits_{x\to +\infty}\widetilde{H}'(x)\cdot \eta,$$
where $\eta \sim N(0,1).$ In other words,
\begin{equation}\sqrt{k_n}\left(H(X_{(n-k_n)}) \frac{S_x(a_{n/k_n},\gamma)}
{S_x(X_{(n-k_n)},\gamma)} - H(a_{n/k_n})\right)  \xrightarrow {d}
 \lim\limits_{x\to +\infty}\frac{\partial\left(\frac{H(x)}{S_x(x, \gamma)}\right)}{\partial x}\cdot N(0,1),\label{morzh}\end{equation}
where we note especially that the function $\frac{H(x)}{S_x(x, \gamma)}$ does not depend on $k_n$ by definition of $H(x).$
Find $\lim\limits_{x\to +\infty}\frac{\partial\left(\frac{H(x)}{S_x(x, \gamma)}\right)}{\partial x}.$ It follows
from the definition of the function $H(x),$ that

$$\widetilde{H}(x) = \frac{H(x)}{S_x(x,\gamma)} = \frac{u}{S_{x\gamma}(x, \gamma)}
\left(\frac{\int\limits_{x}^{\infty}S_{\gamma}(y,\gamma) \exp(-S(y,\gamma)) dy}
{ \int\limits_{x}^{\infty} \exp(-S(y,\gamma))dy} - \right.$$

$$\left. - S_{\gamma}(x,\gamma) - \frac{S_{x\gamma}(x,\gamma)}{S_x(x,\gamma)} -
 \frac{S_{xx\gamma}(x, \gamma)}{(S_x(x, \gamma))^2} + 2 \frac{S_{xx}(x, \gamma)S_{x\gamma}(x, \gamma)}{(S_x(x, \gamma))^3}\right).$$
Denote

$$V(x, \gamma) = \frac{\int\limits_{x}^{\infty}S_{\gamma}(y,\gamma) \exp(-S(y,\gamma)) dy}
{ \int\limits_{x}^{\infty} \exp(-S(y,\gamma))dy} - S_{\gamma}(x,\gamma) - $$
$$- \frac{S_{x\gamma}(x,\gamma)}{S_x(x,\gamma)} - \frac{S_{xx\gamma}(x, \gamma)}
{(S_x(x, \gamma))^2} + 2 \frac{S_{xx}(x, \gamma)S_{x\gamma}(x, \gamma)}{(S_x(x, \gamma))^3}.$$
In the sequel we will omit the arguments of the function $S(x, \gamma),$ $V(x, \gamma)$ and its derivatives,
if they are equal to $x$ and $\gamma$ respectively.
Using (\ref{i2i1}) and the fact that under the conditions of Theorem 2 $\frac{S_{x\gamma}}{S_xS_\gamma}
 = O\left(\frac{S_{xx}}{S_x^2}\right),$ we obtain

$$V = - \frac{S_{xx}S_{x\gamma}}{S_x^3} + o\left(\frac{S^2_{xx}}{S_x^4}\right).$$
So \[\widetilde{H}'(x) = \frac{u V_x}{S_{x\gamma}} - \frac{u V S_{xx\gamma}}{(S_{x\gamma})^2}.\]
Firstly find the asymptotics of the second summand in the last equation
$$\frac{u V S_{xx\gamma}}{(S_{x\gamma})^2} =
- \frac{u S_{xx}}{S^2_x}\frac{S_{xx\gamma}}{S_x S_{x\gamma}} + o\left(\frac{S^2_{xx}}{S_x^4}
\frac{S_{xx\gamma}}{(S_{x\gamma})^2}\right).$$
As it is mentioned previously, under the conditions of Theorem 2 we have
$\lim\limits_{x\to\infty} \frac{S_{xx}}{S^2_x} = 0.$ Then from (\ref{pavian}) imply
 $\lim\limits_{x\to\infty}\frac{S_{xx\gamma}}{S_{x\gamma}S_x} = 0.$
So it is necessary to prove $\lim\limits_{x\to\infty}\frac{S_{xx\gamma}}{(S_{x\gamma})^2} = 0$ to complete the proof of the fact that $\lim\limits_{x\to\infty}\frac{u V S_{xx\gamma}}{(S_{x\gamma})^2} = 0,$
It follows from the L'Hospital rule and the regularity conditions that

$$\lim\limits_{x\to\infty}\frac{S_{xx\gamma}}{(S_{x\gamma})^2} =
\lim\limits_{x\to\infty}\frac{S_{xx\gamma}/S_{x\gamma}}{S_{x\gamma}} =
\lim\limits_{x\to\infty}\frac{\ln S_{x\gamma}}{S_\gamma} \lim\limits_{x\to\infty}
\frac{\ln S_{x\gamma}}{S_\gamma} = 0,$$
that completes the proof of concerned fact.
Now find the asymptotics of the function $\frac{u V_x}{S_{x\gamma}}$ as $x\to+\infty.$
The search of the explicit form of the function $V_x$ requires the complicated calculations, so we omit them and write the answer immediately:
$$\frac{u V_x}{S_{x\gamma}} = - \frac{u S_{xx}}{S_x^2} + o\left(\frac{S_{xx}}{S_x^2}\right),$$
i.e. $\lim\limits_{x\to\infty}\frac{u V_x}{S_{x\gamma}} = 0.$ Hence we derive
 $\lim\limits_{x\to\infty}\widetilde{H}'(x) = 0,$ and from (\ref{morzh}) imply
\begin{equation}\sqrt{k_n}\left(H(X_{(n-k_n)}) \frac{S_x(a_{n/k_n},\gamma)}
{S_x(X_{(n-k_n)},\gamma)} - H(a_{n/k_n})\right)  \xrightarrow {d} 0. \label{los}\end{equation}
From the mentioned previously the continuity theorem for the function $\frac{1}{S_x}$
(where $b_n = \frac{1}{\sqrt{k_n}S_x(a_{n/k_n}, \gamma)}$ and $\eta_n = \sqrt{k_n}S_x(a_{n/k_n},\gamma)(X_{(n-k_n)} - a_{n/k_n}), $ as before)
$$ \frac{\frac{1}{S_x(a_{n/k_n} + b_n\eta_n, \gamma)} - \frac{1}{S_x(a_{n/k_n}, \gamma)}}{b_n}
\xrightarrow{d} \lim\limits_{x\to\infty}\frac{S_{xx}}{S_x^2} \xi,$$
where $\xi\sim N(0,1).$ In other words,
\begin{equation}\sqrt{k_n}\left(\frac{S_x(a_{n/k_n}, \gamma)}{S_x(X_{(n-k_n)}, \gamma)}
 - 1\right) \xrightarrow{d} 0,\label{vepr}\end{equation}
since under the conditions of Theorem 2, $\lim\limits_{x\to\infty}\frac{S_{xx}}{S_x^2}=0.$ From (\ref{i2i1}) imply
 $$H(x) = - \frac{S_{xx}}{S^2_x} + o\left(\frac{S^2_{xx}}{S^4_x} \frac{S_x}{S_{x\gamma}}\right).$$
Since as it is stated previously $\frac{S_x}{S_{x\gamma}} =R(S)$ as $x\to+\infty,$ where $R(S)$ is some slowly varying function, аnd
 $\frac{S_{xx}}{S^2_x} = O(S^{-\delta})$ for some $\delta,$ $0< \delta\leq 1,$ then
  $\lim\limits_{x\to\infty} H(x) = 0.$ Hence using (\ref{vepr}) we obtain
$$\sqrt{k_n} H(X_{(n-k_n)})\left(\frac{S_x(a_{n/k_n}, \gamma)}{S_x(X_{(n-k_n)}, \gamma)} - 1\right) \xrightarrow{d} 0.$$
Combining the previous result and (\ref{los}) we derive finally

$$\sqrt{k_n}\left(H(X_{(n-k_n)}) - H(a_{n/k_n})\right)  \xrightarrow {d} 0,$$
that completes the proof of Theorem 2. $\blacksquare$

\section{References}
\begin{enumerate}
\item [\bf 1]{Antle, C. E., Bain, L. J. (1969). A Property of Maximum Likelihood Estimators of Location and Scale Parameters.
\emph{SIAM Review} \textbf{11 (2)}, 251--253.}

\item [\bf 2]{Billingsley P. (1999). \emph{Convergence of probability measures.} Second edition. Wiley.}

\item [\bf 3] {Borovkov, А.А. (1984). \emph{Mathematical statistics} (in Russian). Мoscow: Nauka.}

\item [\bf 4] {Chibisov, D.M. (2009). \emph{Lectures by asymptotic theory of rank criteria} (in Russian). Lecture courses SEC, 14, MIAS, Мoscow, 3–174.}
\item [\bf 5] {Dey, A. K., Kundu, D. (2012). Discriminating
between the Weibull and log-normal distributions for Type-II
censored data. \emph{Statistics}, \textbf{46 (4)}, 197--214.}

\item [\bf 6] {Dumonceaux, R., Antle, C.E. (1973). Discrimination between the log-normal and the Weibull distributions.
\emph{Technometrics}\textbf{ 15 (4)}, 923--926.}
\item [\bf 7] {Dumonceaux, R., Antle, C.E., Haas, G. (1973). Likelihood ratio test for discrimination between two models
with unknown location and scale parameters. \emph{Technometrics} \textbf{15 (1)},
19--27.}

\item [\bf 8] {Embrechts, P., Kluppelberg, C., Mikosch, T. (1997). \emph{Modelling Extremal Events for Insurance and Finance}. Berlin: Springer.}

\item [\bf 9] {Fedoryuk, М.V. (1977). \emph{The pass method} (in Russian). Мoscow: Nauka.}

\item [\bf 10] {Fereira A., de Haan L. (2006).
\emph{Extreme value theory. An introduction.} Springer Series in
Operations Research and Financial Engineering. N. Y.: Springer.}
\item [\bf 11] {Fraga Alves, I., de Haan, L., Neves, C. (2009). A test procedure for detecting super-heavy tails.
\emph{Journal of Statistical Planning and Inference}, \textbf{138(2),} 213--227.}
\item [\bf 12] {Gardes, L., Girard, S., Guillou, A. (2011). Weibull tail-distributions revisited:
a new look at some tail estimators. \emph{Journal of Statistical Planning and Inference}, \textbf{141(1),} 429--444.}
\item [\bf 13] {Gupta, R.D., Kundu, D. (2003). Discriminating between Weibull and generalized
exponential distributions. \emph{Computational Statistics and Data
Analysis} \textbf{43 (2)}, 179--196.}

\item [\bf 14] {Gupta, R.D., Kundu, D., Manglick, A. (2001). Probability of correct selection of Gamma versus GE or
Weibull versus GE based on likelihood ratio test. Technical
Report, The University of New Brunswick, Saint John.}
\item [\bf 15] {Kundu, D., Raqab, M. Z. (2007). Discriminating Between the Generalized
Rayleigh and Log-Normal Distribution. \emph{Statistics} \textbf{41 (6)}, 505--515.}
\item [\bf 16] {Rodionov, I.V. (2014). Statistical analysis and tests of hypotheses of distributions of time series extrema,
 PhD thesis (in Russian).
 Moscow State University.}
\item [\bf 17] {Roussas, G.G. (1972)\emph{ Contiguity of Probability Measures. Some Applications in Statistics.} Cambridge University Press.}

\end{enumerate}
\end{document}